\documentclass[twoside]{article}
\usepackage[utf8]{inputenc}
\usepackage{amsmath}
\usepackage{amsfonts}
\usepackage{amssymb}
\usepackage{graphicx}
\usepackage{color}
\usepackage[normalem]{ulem}
\usepackage[left=2cm,right=2cm,top=2cm,bottom=2cm]{geometry}
\newenvironment{Proof}{\noindent{\sc Proof.}}{\qed}
\newtheorem{theorem}{Theorem}[section]
\newtheorem{lemma}{Lemma}[section]

\newtheorem{rem}{Remark}[section]
\newtheorem{definition}{Definition}[section]

\newtheorem{uda}{Example}[section]
\newcommand{\qed}{\hfill$\Box$\par\medskip}

\def\bhag#1{\noindent
\setcounter{equation}{0}
\section{#1}
}

\def\RR{{\mathbb R}}
\def\CC{{\mathbb C}}
\def\ZZ{{\mathbb Z}}

\def\PPI{{{\rm I}\kern-1pt\Pi}}
\def\SS{{\mathbb S}}

\def\a{\alpha}
\def\b #1;{{\bf #1}}

\def\e{\epsilon}

\def\O{{\cal O}}

\def\YY{\mathbb{Y}}

\def\derf#1#2{{#1}^{(#2)}}

\def\esssup{\mathop{\hbox{{\rm ess sup}}}}

\def\be{\begin{equation}}
\def\ee{\end{equation}}
\def\bea{\begin{eqnarray}}
\def\eea{\end{eqnarray}}
\def\eref#1{(\ref{#1})}
\def\disp{\displaystyle}

\def\donchitre#1#2{\vskip 6.5cm\noindent
\parbox[t]{1in}{\special{eps:#1.eps x=6.5cm y=5.5cm}}
\hbox to 7cm{}\parbox[t]{0.0cm}{\special{eps:#2.eps x=6.5cm y=5.5cm}}}

\def\tn{|\!|\!|}
\def\XX{{\mathbb X}}
\def\BB{{\mathbb B}}

\def\bs#1{{\boldsymbol{#1}}}

\def\supp{\mathsf{supp\ }}

\title{
A unified framework for  harmonic analysis of functions on directed graphs and changing data }
\author{
H.~N.~Mhaskar\thanks{Department of Mathematics, California Institute of Technology, Pasadena, CA 91125;
Institute of Mathematical Sciences, Claremont Graduate University, Claremont, CA 91711. The research of this author is supported in part by ARO Grant W911NF-15-1-0385.
\textsf{email:} hrushikesh.mhaskar@cgu.edu. }       }
\date{}
\begin{document}
\maketitle
\begin{abstract}
We present a general framework for studying harmonic analysis of functions in the settings of various emerging problems in the theory of diffusion geometry. 
The starting point of the now classical diffusion geometry approach is the construction of a kernel whose discretization leads to an undirected graph structure on an unstructured data set. 
We study the question of constructing such kernels for directed graph structures, and argue that our construction is essentially the only way to do so using discretizations of kernels. We then use our previous theory to develop harmonic analysis based on the singular value decomposition of the resulting non--self--adjoint operators associated with the directed graph. 
Next, we consider the question of how functions defined on one space evolves to another space in the paradigm of changing data sets recently introduced by Coifman and Hirn. While the approach of Coifman and Hirn require that the points on one space should be in a known one--to--one correspondence with the points on the other,  our approach allows the identification of only a subset of landmark points. We introduce a new definition of distance between points on two spaces, construct localized kernels based on the two spaces and certain interaction parameters, and study the evolution of smoothness of a function on one space to its lifting to the other space via the landmarks. We develop novel mathematical tools that enable us to study these seemingly different problems in a unified manner.
\end{abstract}

\bhag{Introduction}\label{intsect}
There are many approaches developed during the last decade or so in order to analyze large, unstructured, possibly high dimensional data. The basic idea is to embed the high dimensional data on a low dimensional sub--manifold of the ambient Euclidean space.
The main theme of the research is then to understand the data geometry in terms of the geometric properties of this sub--manifold. Well known techniques in this direction, dimensionality reduction in particular, are  Isomaps \cite{tenenbaum2000global}, maximum variance unfolding (MVU) (also called
 semidefinite programming (SDP)) \cite{weinberger2005nonlinear}, locally linear embedding
 (LLE) \cite{roweis2000nonlinear}, local tangent space alignment method (LTSA) \cite{zhang2004principal}, Laplacian eigenmaps (Leigs) \cite{belkin2003laplacian}, Hessian locally
 linear embedding (HLLE) \cite{david2003hessian},  diffusion maps (Dmaps) \cite{coifmanlafondiffusion}, and randomized anisotropic transform \cite{chuiwang2010}. A recent survey of these methods is given by Chui and Wang in \cite{chuidimred2015}.  An excellent introduction to the subject of diffusion geometry can be found in the special issue \cite{achaspissue} of Applied and Computational Harmonic Analysis, 2006. The application areas are too numerous to mention exhaustively. They include, for example, document analysis \cite{coifmanmauro2006}, 
 face recognition \cite{niyogiface, ageface2011, chuiwang2010}, hyperspectral imaging \cite{chuihyper},  semi-supervised learning \cite{niyogi1, niyogi2}, image processing \cite{donoho2005image, arjuna1},  cataloguing of galaxies \cite{donoho2002multiscale},  and social networking \cite{bertozzicommunity}.
 
The starting point in diffusion geometry is the point cloud.  A point cloud is a set $\mathcal{P}$ of points $\{x_i\}$ in a Euclidean space, together with a
similarity relation $W$, viewed as the matrix of edge weights in  an undirected graph. In the absence of any known structure on the set of points, the similarity relation is constructed by the so--called diffusion matrix. For example, a standard construction for $W$ is given by
$$
W_{i,j}=\exp(-\|x_i-x_j\|^2/\e),
$$
where $\{x_i\}$ is the set of points, 
$\|\cdot\|$ is the Euclidean norm, and $\e$ is a judiciously chosen variance parameter.  The first few eigenvalues and eigenvectors of this matrix (or some related matrix, such as $\mathsf{diag}(\mbox{row sums of }W)-W$) provide the desired low dimensional embedding. Usually, these vectors codify certain identifiable features of the data set. For example, in \cite{chuiwang2010}, the point cloud consists of thumbnail images of the faces of the same person in different orientations. The components of the first non--trivial eigenvector orders the collection according to the angle of rotation. Many such examples can be found in the literature; indeed, the well known concept of kernel PCA is based on this fact. It is proved in \cite{belkin2003laplacian, belkinfound, lafon, singer} that as the data becomes dense, the so--called graph Laplacian based on a judiciously chosen (weighted) adjacency matrix, such as $W$, converges to the Laplace--Beltrami operator on the \textbf{unknown} manifold from which the data is sampled, and likewise for the corresponding eigenvalues/eigenfunctions. 
A deeper insight into the intimate connection between the eigenfunctions and the manifold is explained in \cite{jones2008parameter, jones2010universal}, where it is shown that some of the eigenfunctions define a local coordinate system on the unknown manifold so that the Euclidean distance between the points represented with these coordinates is proportional to the geodesic distance between the points.

The applications of this theory to semi--supervised learning can be formulated as problems of function extension, going beyond the important and difficult question of understanding the data geometry. For example, in semi--supervised learning for classification, the class labels are known only on a small \emph{training subset} $\mathcal{C}\subset \mathcal{P}$, and the objective is to find the class labels for all points in $\mathcal{P}$. The class label for a point $x$ can be viewed as the value of a \emph{target function} $f$ at $x$, and the question is then to extend $f$ from $\mathcal{C}$ to $\mathcal{P}$. Formulated in this manner, one can even think of extending $f$ to the entire manifold $\XX$, including those points which are not in the original data set $\mathcal{P}$. 
Clearly, all the regression problems in learning theory are also problems of function extension/approximation, for example \cite{girosi1990networks}. Indeed, one of the main reasons for the popularity of neural and radial basis function networks in learning theory is their universal approximation property. An important problem in this paradigm is to estimate the fidelity of the approximation scheme at unseen data points. Within the context of diffusion geometry, a rigorous analytic theory to address this question is developed in \cite{mauropap, frankbern, modlpmz, heatkernframe, eignet}, based on an abstract framework formulated in \cite{fasttour}. The constructive algorithm that emerges from this theory is tested in the proof--of--concept experiments on recognition of hand--written digits \cite{mauropap}, the Cleveland heart disease data set, the Wisconsin breast cancer data set, and a new experiment at NIH to predict automatically the local classifications of drusen in patients with age related macular degeneration (AMD) based on multi--spectral fundus images of the retina \cite{compbio}.  

The goal of the present paper is to develop a general framework in which function extension/approximation problems can be studied in the context of two of the newly emerging paradigms in the theory of diffusion geometry, 
namely, the problem of changing data sets, and the problem of analyzing functions on a directed graph.

The problem of analyzing changing data arises in  many applications, for example, modeling of social networks, where the relationships between the people involved may change over time, analysis of financial markets, evolutionary biological questions, and  the analysis of medical tests in patients as they change over a time period in the life of the patient.  In classical diffusion geometry, the eigenfunctions and eigenvalues are based on an existing data set. The out--of--sample extension problem is the problem of incorporating new data as they emerge. This problem can be viewed clearly as the problem of changing data as well.

Our technical motivation for the analysis of changing data sets is the paper \cite{coifmanhirn} of Coifman and Hirn. In this paper, the authors define a diffusion distance between points lying on two different manifolds, as well as a distance between changing manifolds, based on these. One 
important feature of this framework 
is that  the points in the point clouds (respectively, the underlying manifolds) remain the same over different manifolds, and can be identified uniquely in a one--to--one correspondence as the manifold changes. The authors illustrate an application of this framework  in the context of a sequence of hyperspectral images of essentially the same area over time. Each image is treated as a point cloud, where the ``points'' are the spectral curves. Each point in one image is identified with exactly one point in another image at the same  spatial pixel location.

In practical applications, it is often not possible to identify each point on one manifold uniquely with the corresponding point on another manifold. For example, in the analysis of brain MRI images of Alzheimer patients taken over time, the images are taken at  random points of the brain \cite{kim2014multi}. 
Therefore, even if the brain surfaces are the ``same'' in the sense of belonging to the same patient, it is not possible to identify one point in the image uniquely
 with another point on another image. In fact, there are different  number of points in different images.  One may, however, consider the cortical thickness of the brain as a function on an imaginary sphere, approximate its values at certain standard points, such as dyadic points, which can then be identified from one sphere to the other in the sequence of images as they evolve in time. 
 
 In a different example, if one considers each frame of a video to be a point cloud/manifold as is often done in image processing applications of diffusion geometry, each point on one video frame can get mapped to different points on the next frame with different probabilities \cite{videomoeslund2012}. However, it may be reasonable to assume that certain ``landmark'' points can be identified from manifold to manifold as they change.
A related direction of research is manifold matching, where one needs to match features of one manifold with the other to determine whether the two manifolds arise from different representations of the same data. 

The framework developed in this paper enables us to study the lifting of functions from one manifold to another via a set of landmarks rather than requiring a complete one--to--one correspondence between points of the different manifolds.
 
Another important emerging direction is to extend the diffusion geometry paradigm to directed graphs rather than undirected graphs defined by a symmetric matrix such as $W$ above. For example,
graphs arising from many applications such as transportation problems in which traffic flow is restricted to one direction must be directed graphs, in order to allow analysis of traffic along one way streets. In analyzing data about world--wide web links, the direction of the links contains crucial information.  Similarly, in an understanding of blood circulation, the direction of blood flow is important. Other applications of directed graphs include one--way communication problems,  social interactions between people where the relationship is not necessarily symmetric, and scheduling problems.   
 
Directed graphs are, of course, studied in great detail over more than a century. In particular, the notion of symmetric graph Laplacian for a directed graph whose weighted adjancancy matrix is given by a Markov matrix is developed by F. Chung in \cite{chung_directed_laplacian}. A recent well written Ph. d. thesis by Gidelew \cite{gidelew2014topics} develops harmonic analysis based on this Laplacian, and cites several other papers in this direction.
 While this line of research assumes that edge  weights in the graph are already known, it is a separate question to impose a directed graph structure for an unstructured data, analogous to the graph structure implied on undirected data by matrices such as $W$. One may view the asymetric affinity matrices  in \cite{coifman2013bi} as kernels defining a directed graph in the abstract sense. To the best of our knowledge, the only deliberate effort in this direction is the recent paper  \cite{mousazadeh2015embedding} by Mousazadeh and Cohen, where a modified Gaussian kernel is proposed, and an analogue of the limiting Laplace--Beltrami operator is presented. The theory of non--self--adjoint integral and differential operators is quite well developed, e.g., \cite{gohbergkrein, halmosbdintop, marockner1992}. It is therefore natural to obtain the desired graph structure by discretizing a suitable  kernel, deemed most appropriate for the application in question.

Unlike \cite{belkin2003laplacian, belkinfound, lafon, singer} in the case of undirected graphs or \cite{mousazadeh2015embedding} in the case of directed graphs, it is not our intention to focus on specific constructions of the kernels which give rise to some concrete differential operator in the limiting sense. The book \cite{marockner1992} of Ma and R\"ockner gives a very readable account of the interconnection between non--self--adjoint differential operators, the associated Dirichlet forms, Green's kernels, and semi--groups of operators generated by the differential operators. 
Several examples are also discussed in \cite{marockner1992}, and may be adapted for this purpose.
Instead, we wish to enumerate certain properties that such kernels ought to have in order to develop a 
reasonably rich theory of harmonic analysis of functions defined on changing data or directed graphs. We will, however, make some comments on what we consider likely to be a good construction of kernels in these cases, leaving the details in each case to future research. In particular, in the case of construction of a  directed graph structure on data, we will argue that our suggestions provide essentially the only way to have such a structure, as well as enable the construction of such structure to allow meaningful harmonic analysis based on this structure.

To describe our results more precisely, it is convenient to review the setting in the ``classical'' diffusion geometry theory \cite{coifmanmauro2006}. Let $\XX$ be a metric measure space with a metric $d$ and a sigma--finite Borel measure $\mu^*$. 
Let $\{\phi_k\}$  be a complete orthonormal system in $L^2(\mu^*)$, considered  as the eigenfunctions of a Laplician operator defined  by a formally defined kernel
$$
-\sum_{k=0}^\infty \lambda_k^2\phi_k(x)\overline{\phi_k(y)}, \qquad x, y\in\XX,
$$
where $\lambda_k$ is a sequence of non--negative numbers with $\lambda_k\uparrow\infty$ as $k\to \infty$.
For the analysis of a measurable function $f :\XX\to\CC$, Coifman and Maggioni \cite{coifmanmauro2006} propose taking dyadic powers of the ``diffusion'' operator
$$
f\mapsto\int_\XX f(y)\left\{\sum_{k=0}^\infty \exp(-\lambda_k^2t)\phi_k(\cdot)\overline{\phi_k(y)}\right\}d\mu^*(y), \qquad x\in \XX,
$$
for a small value of $t$.
In practice this is done as dyadic powers of a sparse matrix using fast algorithms. Theoretically, the range of the $2^{-2n}$--th power of this operator is approximately equal to the space
$$
\Pi_{2^n\e}=\mathsf{span}\{\phi_k : \lambda_k<2^n\e\}.
$$
for a suitable $\e$.
The corresponding projection operator is given by
$$f\mapsto \int_\XX f(y) \left\{\sum_{k : \lambda_k <2^n\e} \phi_k(\cdot) \overline{\phi_k(y)}\right\}d\mu^*(y). 
$$

For the detection of certain anomalies in the target function, it is desirable to approximate the function in the uniform norm rather than the $L^2(\mu^*)$ norm.
 It is well known that the projection operators are usually unbounded if uniform approximation is desired, even in the very classical case of trigonometric Fourier series.
  In \cite{fasttour, mauropap}, we have initiated a study of function approximation and harmonic analysis in $L^p(\mu^*)$ for values of $p\not=2$; with a special emphasis on uniform approximation.  The important starting point of our theory is therefore to construct localized operators of the form
\be\label{intsigmaopdef}
f\mapsto \int_\XX f(y)\left\{\sum_{k=0}^\infty h\left(\frac{\lambda_k}{n}\right)\phi_k(\cdot)\overline{\phi_k(y)}\right\}d\mu^*(y)
\ee
for a suitable (low pass or band pass) ``filter'' $h$. 

Taking a sufficiently smooth function $h$ as the filter ensures that these operators  are uniformly bounded in every $L^p(\mu^*)$, $1\le p\le\infty$. Using these operators, one can obtain  wavelet--like dyadic representations in all of these spaces, where the behavior of the terms characterizes completely the smoothness properties of the target function $f$. One fascinating feature of this construction is that even if the operator itself is defined using the spectral information about the target function $f$, the localization estimate leads one to the analogues of the well known theorems of Jaffard \cite[Theorems~9.2.1, 9.2.2]{daubbook}; i.e., the local smoothness of the target function
at a point can be characterized completely in terms of the behavior of the terms of the wavelet--like representation in a neighborhood of the point in question. It is proved in \cite{chuiinterp} that such localized harmonic analysis can always be done when appropriate localized kernels are available.

A feature common to both the new directions of research mentioned above that distinguishes the new problems from the classical one is that one has to deal with two  systems jointly. Thus, in the case of changing data or manifold matching, one has two admissible systems $\Xi_1=(\XX_1,d_1,\mu_1^*,\{\lambda_{1,k}\}, \{\phi_{1,k}\})$ and 
$\Xi_2=(\XX_2,d_2,\mu_2^*,\{\lambda_{2,k}\}, \{\phi_{2,k}\})$ (cf. Definition~\ref{systemdef} for a precise definition). The objective of interest in \cite{coifmanhirn} was to define a ``diffusion distance'' between points on the two manifolds using the respective heat kernels on the two manifolds. In our paper, the goal is to study how a function defined on one system evolves into another function on the other system via an interaction between the various parameters defining the two systems. We will assume the interaction only in an abstract sense, without requiring that each point on one manifold should be identified uniquely with another
point on the other manifold. Moreover, the construction and theory is general enough to deal with situations where a point on one manifold may go to different points on the other manifold according to some probability distribution.

The case of directed graphs is much easier in some sense. A directed graph is generated in this paradigm by the discretization of a non--self--adjoint kernel operator.  Thus, we again have two admissible systems, where the underlying  measure spaces, measures, and distances are the same in both the systems. In analogy to the classical diffusion setting, it is convenient to take this kernel to be a suitable kernel related to some differential operator on the underlying space, and the sequence $\{\lambda_k\}$ to be the singular values of the differential operator. The two systems arise by taking for $\{\phi_{1,k}\}$, $\{\phi_{2,k}\}$ the left and right singular functions of the operator. With this insight, there is essentially only one way to construct such operators, stemming from the polar decomposition of the non--self--adjoint operator involved, and the theory which we have already developed can be applied directly with minor changes.

The basic objective of this paper is to obtain the analogues of the localized operators analogous to those in \eref{intsigmaopdef} for harmonic analysis in these settings. Obviously, it cannot hold that these operators  applied to a function on say $\XX_2$, will converge to the function itself -- the result of this application is not even a function on $\XX_2$. However, under certain conditions, one can obtain a lifting operator from one space to the other so as to result
in a lifting of the original starting function to a corresponding function on the other space, and the smoothness properties of this lifted function can be described using those of the starting function. A motivating model here is 
the classical theorems in the theory of partial differential equations (or integral equations), such as trace theorems, inverse theorems, or embedding theorems \cite{nikolskii}.

The important technical tool for achieving the unified framework is 
to prove what is known as a  Tauberian theorem. To describe this idea, we make an observation in the classical setting.  If one considers the univariate function defined by
$$
\mu(u)=\mu(x,y;u)=\sum_{k:\lambda_k <u}
\phi_k(x)\overline{\phi_k(y)}, \qquad u>0,
$$
then $\mu$ has bounded total variation on compact sub--intervals of $[0,\infty)$, and the projection kernel
 can be expressed formally in the form of the Stieltjes integral
$$
\sum_{k:\lambda_k<n}\phi_k(x)\overline{\phi_k(y)}=\int_0^\infty h(u/n)d\mu(u),
$$
where $h(t)=1$ if $0\le t<1$, and $=0$ otherwise. The kernel of the operators defined in \eref{intsigmaopdef} can also be described in the same way with different choices of $h$.
Likewise, the heat kernel can be expressed formally as a Stieltjes integral
$$
\int_0^\infty \exp(-u^2t)d\mu(u).
$$
Our theorems in \cite{mauropap, frankbern} on the  localization properties of the kernel of the operators defined in \eref{intsigmaopdef} can therefore be viewed as  statements about the properties of one integral transform of $\mu$ in terms of those of another transform.  Such a theorem is called a Tauberian theorem; a famous example being the Wiener--Ikehara Tauberian theorem regarding the density of translates of a fixed function in $L^1(\RR)$. The history of such theorems in the  settings of classical harmonic analysis spans over a century, and is summarized excellently in a recent book \cite{korevaar_taub_book} by Korevaar. Early versions of the Tauberian theorems of the form which we will prove in this paper were obtained by Bochner and Chandrasekharn \cite{bochner_chadra_book}. 

Viewed in this abstraction, the specific nature of $\mu$ becomes irrelevant, in particular, the dependence of $\mu$ on $x, y$ or systems $\{\lambda_k\}$, $\{\phi_k\}$ does not play any role anymore. Therefore, the problem of constructing localized kernels in two or more systems reduces to the problem of proving an appropriate Tauberian theorem. While the  classical development of such theorems is focused on the limiting or asymptotic behavior of a measure $\mu$ itself near $\infty$ or some summability operator applied to $\mu$, our theorem has a very different character. We are interested in obtaining some very specific bounds on one transform when some other very specific bounds on the other transform are known. In particular, the proof relies heavily upon both complex variable theory as well as distribution theory.

We summarize the main new contributions of this paper as follows:
\begin{itemize}
\item We point out a very simple way to construct kernels to embed a directed graph into a manifold, and argue that this is essentially the \textbf{only way} to construct such kernels.
\item We develop a theory for harmonic analysis of functions on directed graphs, whether using existing kernels as the starting point, or by constructing specially designed kernels.  
\item We develop a theory of lifting a function on one data space to another using only the identification of certain interaction parameters connecting the two  spaces.
\item We develop a mathematical theory that enables us to treat both of these seemingly unrelated problems in a unified manner. Other contemplated applications of this theory include an analysis of solutions of partial differential equations with non--self--adjoint operators, bit representation of signals on Euclidean spaces, and smoothness--preserving image completion.
\end{itemize}

In Section~\ref{dirgraphsect}, we will present our suggestions for the construction of kernels whose discretizations yield the desired directed graph structure on an unstructured data. We will argue that our construction is essentially the only one based on the idea of discretizing operators, barring some esoteric counterexamples. We will then discuss harmonic analysis
on the systems generated by these non--self--adjoint kernels. The results are partly a review of our previous results, and partly a bridge between these and the results on changing data, which will be described in Section~\ref{twosystemsect}. The proofs of all the results, including the Tauberian theorem will be given in Section~\ref{pfsect}.

\bhag{Directed graphs}\label{dirgraphsect}

In this section, we focus on a discussion of harmonic analysis on directed graphs. In Sub--section~\ref{digraphconstsect}, we discuss the construction of the systems themselves, either from a pre--existing structure, or from the point of view of constructing kernels in order to embed the given data as a directed graph on points on an ambient manifold. We will point out that there is essentially only one way to make such constructions. In Sub--section~\ref{digraphharmanalsect}, we assume that the system is constructed, and study the harmonic analysis of functions on the graphs based on the system assuming certain conditions. The results are obtained as immediate consequences of our results in \cite{mauropap, frankbern}, and serve both as a review of these results and a prelude to the results to be presented in Section~\ref{twosystemsect}. We note that in the continuous setting, it is natural to think of a directed graph as a non--self--adjoint operator. 
Therefore, it is not necessary that this operator must map a Hilbert space to itself. For example, in applications in photo--acoustic tomography one 
encounters operators defined on one Hilbert space, taking values in a different Hilbert space \cite{filbirmadychacha}. 
We will explain in Remark~\ref{diffhilbertrem} how the analysis presented in this section can be extended easily to such situations.  However, for simplicity, we will restrict ourselves in this section to the case of one space.

\subsection{Construction of the graph structure}\label{digraphconstsect}

To motivate our discussion, let us recall that a directed graph with the vertex set $\{x_i\}_{i=1}^M$, and edge weights $w(x_i,x_j)$  is an $M\times M$ non--symmetric matrix $W$ whose entries are $w(x_i,x_j)$. Then the polar form for $W$ \cite[Corollary~7.3.3]{horn_johnson_book} is $PU$, where $P$ is a uniquely determined positive semi--definite matrix with the same rank as $W$, and $U$ is unitary and real. If $W$ has full rank, then $U$ is determined uniquely. In this case, efficient algorithms to compute the polar decomposition are given recently by Nakatsukasa, Bai, and Gygi \cite{nakatsukasa2010optimizing}.

Writing $P$ in the form
$$
P(x_i,x_j)=\sum_{k=0}^{M-1}\lambda_k\phi_k(x_i)\phi_k(x_j),
$$
for some orthogonal vectors $\phi_k$, the singular value decomposition of $W$ is given by
$$
W(x_i,x_j)=\sum_{k=0}^{M-1}\lambda_k\phi_k(x_i)(U^*\phi_k)(x_j).
$$
Thus, if $f$ is any function on the points $\{x_i\}$, then 
$$
(Wf)(x_i)=(P(Uf))(x_i).
$$
Therefore, the harmonic analysis of $f$ based on the singular vectors of $W$ is the same as the harmonic analysis of $Uf$ based on the eigenvectors of $P$, as in \cite{fasttour, hammond}. Moreover, this is essentially the only shape such a harmonic analysis can take.

We digress to make a comparison with the  construction of the graph Laplacian for directed graphs proposed in \cite{chung_directed_laplacian}, where the role of $W$ is played by a non--symmetric matrix denoted in \cite{chung_directed_laplacian} by $\Phi^{1/2}P\Phi^{-1/2}$. Writing in this discussion, $I$ for the identity matrix, the graph Laplacian in \cite{chung_directed_laplacian} is then given by 
$$
I-\frac{W+W^*}{2}=2\left(I-\left(\frac{W+I}{2}\right)\left(\frac{(W+I)^*}{2}\right)\right) -\frac{1}{2}(I-WW^*),
$$ 
which is the same for  the graph Laplacian corresponding to the underlying undirected graph given by $(W+W^*)/2$. The quantity $P^2=WW^*$ in our notation refers to a different undirected graph. If $W_{i,j}$ denotes the transition probability  of going from state $i$ to state $j$, $(WW^*)_{i,j}$ denotes the probability of nodes $i$ and $j$ leading to some common state.
The matrix $(1/2)(W+I)((1/2)(W+I))^*$ can also be given a similar interpretation, where self--loops are added to the original graph at each vertex.

In the case of infinite graphs, the same analysis and logic works, provided the infinite matrix of edge weights represents a compact operator. The operator $U$ in this case may only be a partial isometry rather than a unitary operator, but this detail does not affect our analysis in this section. In this case, we can view the matrix as the kernel, and the underlying measure to be an atomic measure \cite{halmosbdintop}.
If we do not have a choice of the weight matrix, then clearly, the  problem of finding  the polar decomposition is equivalent to the problem of finding the singular value decomposition. However, the basic fact that the resulting harmonic analysis will be the same as the harmonic analysis of $Uf$ based on the symmetric component is still the same.

In the case when there is no preconceived notion of what the weight matrix ought to be, and we have to impose the structure as in the case of classical diffusion geometry, these observations make the problem of kernel construction almost trivial. 
We start with any kernel based, self--adjoint, positive semi--definite operator, with any familiar kernel used in classical diffusion geometry that leads to a compact operator, such as the Gaussian, or some Green's kernel of a desired elliptic differential operator, take a composition with an appropriate partial isometry (even a unitary operator for simplicity), and discretize the resulting operator. 
This approach leads to a very flexible construction, where also the unitary operator may be treated as a parameter to be chosen according to the application. In particular, if $P$ is an operator that commutes with a self--adjoint  operator $L$  then it is obvious that $W$ and $LU$ commute in the sense that $(LU)^*W=W^*(LU)$ and $W(LU)^*=(LU)W^*$.

We observe that this is essentially the \textbf{only} way to construct a reasonable directed graph structure in the paradigm of diffusion geometry; i.e., by discretizing kernel based operators. For example, if the non--self--adjoint operator is a Hilbert--Schmidt operator, then a polar decomposition is guaranteed \cite{gohbergkrein}, and the fact that the self--adjoint operator is also Hilbert--Schmidt, and hence, kernel based, is observed, for example, in \cite{halmosbdintop}.

Thus,  a simple example of the construction in \cite{mousazadeh2015embedding} is the following. We wish to construct an operator on a Euclidean space $\RR^q$. Fixing $z^*\in\RR^q$ and $t>0$, we define
$$
(W_tf)(x)=(2\pi t)^{-1/2}\int_{\RR^q} \exp\left(-\frac{\|x-y-z^*\|^2}{t}\right)f(y)dy=(2\pi t)^{-1/2}\int_{\RR^q} \exp\left(-\frac{\|x-y\|^2}{t}\right)f(y-z^*)dy.
$$
It is readily seen that this has the polar decomposition as above where the positive definite self--adjoint operator is defined by the Gaussian kernel, and the unitary operator is defined by $f\mapsto f(\cdot-z^*)$. Greater flexibility is obtained by treating $z^*$ as a parameter in addition to $t$. In 
\cite{mousazadeh2015embedding}, 
the integral is taken over an unknown manifold,  the inner product $z^*\cdot (x-y)$ is replaced by a vector field operating on $x-y$, and some numerical experiments are presented to illustrate the idea.

\subsection{Harmonic analysis}\label{digraphharmanalsect}
In this section, we will assume that the necessary structures are available, and develop harmonic analysis based on these. This section may be viewed partly as a review of our previous work on one system, and partly as an example of the results in Section~\ref{twosystemsect}. For this reason, some of the definitions will be more general than necessary in the context of directed graphs.

First, we describe some terminolgy.

If $\XX$ is a measure space and $\nu$ is a (signed complex valued, or positive) measure on $\XX$, then for measurable $f:\XX\to\CC$, we define
\be\label{lpnormdef}
\|f\|_{\XX,\nu,p}=\left\{\begin{array}{ll}
\disp\left\{\int_\XX |f(x)|^pd|\nu|(x)\right\}^{1/p}, &\mbox{ if $1\le p<\infty$},\\[2ex]
|\nu|-\esssup_{x\in\XX}|f(x)|, & \mbox{ if $p=\infty$}.
\end{array}\right.
\ee 
The space $L^p(\XX,\nu)$ consists of all (equivalence classes of) functions $f :\XX\to\CC$ for which $\|f\|_{\XX,\nu,p}<\infty$. If the space $\XX$ is clear from the context, then we will drop its mention from the notation, and similarly if the measure $\nu$ is clear from the context, we will drop its mention as well, as long as the mention of $\XX$ is dropped.
Thus, for example, $L^p(\mu^*)$ refers to $L^p(\XX,\mu^*)$.

We recall that if $\XX$ is a measure space with a positive measure $\mu^*$, a system $\{\phi_k\}_{k=0}^\infty\subset L^2(\mu^*)$ is called a \textbf{Bessel system} if there is a  subset $\mathcal{D}$ of  $L^1(\mu^*)\cap L^\infty(\mu^*)$, dense in $L^2(\mu^*)$, and a functional $\mathcal{N}$ on $L^2(\mu^*)$ such that
\be\label{besselsystdef}
\sum_{k=0}^\infty \left|\int_\XX f(x)\overline{\phi_k(x)}d\mu^*(x)\right|^2 \le \mathcal{N}(f), \qquad f\in \mathcal{D}.
\ee
 For example, if $\XX$ is a complete Riemannian manifold without boundary, and $\mu^*$ is its volume measure, $F$ is a vector field defined on $\XX$, and $\{\phi_k\}$ are the eigenfunctions of the Laplace-Beltrami operator on $\XX$ then the Green's formula implies that $\{F\phi_k\}$ is a Bessel system. 

\begin{definition}\label{systemdef}
Let $\XX$ be a locally compact quasi--metric space with a sigma--finite positive Borel measure $\mu^*$, $d$ be the quasi--metric on $\XX$, $\{\lambda_k\}$ be a non--decreasing sequence of positive numbers, with $\lambda_0=0$ and $\lim_{k\to\infty}\lambda_k=\infty$. Let $\{\phi_k\}$ be a Bessel system of functions in $L^2(\XX,\mu^*)$, such that each $\phi_k$ is continuous, bounded, and integrable on $\XX$. Then $\Xi=(\XX, d, \mu^*, \{\lambda_k\}, \{\phi_k\})$ will be called an \textbf{admissible system}.
\end{definition} 

In this section, we assume that a directed graph is represented in abstract by an admissible system\\
 $\Xi=(\XX, d, \mu^*, \{\lambda_k\}, \{\phi_k\})$ (cf. Definition~\ref{systemdef}), where $\{\phi_k\}$ is an orthonormal set of functions in $L^2(\mu^*)$, together with a partial isometry operator $U$ on $L^2(\mu^*)$. The measure $\mu^*$ may be a discrete measure, and we do not assume that the system is constructed with any particular choice of kernels. Therefore, our theorems are applicable in the settings of finite or infinite pre--defined directed graphs or in the setting where a directed graph structure is constructed in order to embed the graph on a manifold, as long as the required assumptions to be enumerated in Theorem~\ref{digraphanaltheo}  are satisfied. The system $\{\psi_k=U^*\phi_k\}$ is necessarily also an orthonormal system, so that the system $\Xi'=(\XX, d, \mu^*, \{\lambda_k\}, \{\psi_k\})$ is also an admissible system. 
 As we pointed out in Section~\ref{digraphconstsect}, the harmonic analysis of a function $f$ on $\XX$ in the setting of this abstract directed graph is actually a decomposition of $Uf$ into a convergent series of components band limited to different frequency bands with respect to the system $\{\phi_k\}$, based on the data on $f$ with respect to the system $\{\psi_k\}$. We will show that the behavior of these components characterize the  smoothness properties of $Uf$.

In the absence of any differentiability structure, one needs to define the notion of smoothness carefully. We will do this for the system $\Xi$, but will need to use it also for the system $\Xi'$. 
First, we define for $f\in L^1(\mu^*)+L^\infty(\mu^*)$,
\be\label{genfourcoeffdef}
\hat{f}(k)=\hat{f}(\Xi; k)=\int_\XX f(y)\overline{\phi_k(y)}d\mu^*(y), \qquad k=0,1,\cdots.
\ee
The analogue of the classical Sobolev classes with smoothness index $\gamma$ is the class of all functions in $L^2(\mu^*)$ for which
$$
\sum_{k=0}^\infty (k^2+1)^{\gamma}|\hat{f}(k)|^2 <\infty.
$$
This definition has the advantage that it does not depend upon any differentiability structure on $\XX$. However, even in the most classical case of the trigonometric Fourier analysis, it is well known that the Fourier coefficients do not characterize the smoothness of the function in spaces other than the space of $2\pi$--periodic functions square integrable on $[-\pi,\pi]$. So, this definition does not quite work in the case of other $L^p$ spaces. To modify this definition, we will use the notion of the degree of approximation, $E_{n,p}(\Xi; f)$ defined in \eref{degapproxdef} below.

For $n>0$ (not necessarily integer), let
\be\label{polyspacedef}
\Pi_n=\Pi_n(\Xi)=\mathsf{span}\{\phi_k : \lambda_k<n\},
\ee
\be\label{degapproxdef}
E_{n,p}(f)=E_{n,p}(\Xi; f)=\inf_{P\in\Pi_n}\|f-P\|_{\XX,\mu^*,p}, \qquad f\in L^p(\XX,\mu^*), \ 1\le p\le \infty.
\ee
The class of all $f$ for which $E_{n,p}(f)\to 0$ as $n\to\infty$ will be denoted by $X^p=X^p(\Xi)$.
In view of Parseval identity we observe that if $f\in X^2(\Xi)$, then
$$
E_{n,2}^2(f)=\sum_{k: \lambda_k\ge n}|\hat{f}(k)|^2.
$$
Hence, using a routine dyadic sum argument, it is easy to deduce that
$$
\sum_{k=0}^\infty (k^2+1)^{\gamma}|\hat{f}(k)|^2 \sim \sum_{j=0}^\infty 2^{2j\gamma}E_{2^j,2}(f)^2.
$$
(Here, and in the sequel, we use the following convention regarding constants. The symbols $c, c_1,\cdots$ will denote generic positive constants independent of the obvious variables, such as $f$ in the above equation. The values of these constants may be different at different occurrences, even within the same formula. The symbol $A\sim B$ will mean that $c_1A\le B\le c_2B$.)
This observation prompts the following definition.
For $0<\rho\le\infty$ and $0<\gamma<\infty$, the Besov space $B_{p,\rho,\gamma}=B_{p,\rho,\gamma}(\Xi)$ is the class of all $f\in X^p$ such that the quantity $\|f\|_{p,\rho,\gamma}<\infty$, where
\be\label{besovspacedef}
\|f\|_{p,\rho,\gamma}=\|f\|_{\Xi, p, \rho,\gamma}=\left\{\begin{array}{ll}
\disp\|f\|_{\XX,\mu^*,p}+\left(\sum_{j=0}^\infty 2^{\gamma\rho j}E_{2^j,p}(f)^\rho\right)^{1/\rho}, &\mbox{ if $0<\rho<\infty$},\\[1ex]
\|f\|_{\XX,\mu^*,p}+\sup_{j\ge 0} 2^{\gamma j}E_{2^j,p}(f), &\mbox{ if $\rho=\infty$}.
\end{array}\right.
\ee
The index $\gamma$ may be thought of as the smoothness parameter, the parameter $\rho$ adds an additional level of refinement to the smoothness classes. 
  It is well known in approximation theory that the smoothness of a function $f\in X^p$ as measured by $K$--functionals can be characterized in terms of $E_{n,p}(f)$. In the context of diffusion geometry, we have given such theorems in \cite{mauropap, frankbern}. 
  
It is convenient to define a sequence space $\mathsf{b}_{\rho,\gamma}$ consisting of all sequences $\mathbf{a}=\{a_k\}_{k=0}^\infty$ for which
\be\label{seqspacedef}
\|\mathbf{a}\|_{\mathsf{b}_{\rho,\gamma}}=\left\{\begin{array}{ll}
\disp\left(\sum_{j=0}^\infty 2^{\gamma \a j}|a_j|^\rho\right)^{1/\rho}, &\mbox{ if $0<\rho<\infty$,}\\[1ex]
\sup_{j\ge 0}2^{j\gamma}|a_j|, &\mbox{ if $\rho=\infty$},
\end{array}\right.
\ee
is finite. Thus, $B_{p,\rho,\gamma}$ consists of $f\in X^p$ such that $\{E_{2^j,p}(f)\}_{j=0}^\infty \in \mathsf{b}_{\rho,\gamma}$.

Our analysis of functions on $\XX$
is based on spectral information of the form $\{\hat{f}(\Xi';k)\}_{k=0}^\infty$, but the smoothness of $Uf$ is measured in terms of the Besov spaces $B_{p,\rho,\gamma}(\Xi)$. 

To obtain the wavelet--like representation, we introduce some further terminology.

For a function $H :[0,\infty)\to [0,\infty)$, we define formally a kernel 
\be\label{digraphkerndef}
\Phi_n(\Xi,\Xi';H,x,y)=\sum_{k=0}^\infty H\left(\frac{\lambda_k}{n}\right)\phi_k(x)\overline{\psi_k(y)}, \qquad x, y\in \XX.
\ee

\begin{definition}\label{lowpassdef}
A function $h:\RR\to [0,1]$ is called a \textbf{low pass filter} if $h$ is even, $h(u)=1$ if $|u|<1/2$, $h(u)=0$ if $|u|\ge 1$, and $h$ is non--increasing on $[0,\infty)$. 
\end{definition}

With a low pass filter $h$, we next define  a reconstruction operator for $f\in L^1(\mu^*)+L^\infty(\mu^*)$, $x\in\XX$ by
\be\label{digraphsigmadef}
\sigma_n(\Xi,\Xi';h,f,x)=\sum_{k=0}^\infty h\left(\frac{\lambda_k}{n}\right)\hat{f}(\Xi';k)\phi_k(x) =\int_\XX \Phi_n(\Xi,\Xi';h,x,y)f(y)d\mu^*(y), 
\ee
and the analysis operators
\be\label{digraphtauopdef}
\tau_j(\Xi,\Xi';h,f,x)=\left\{\begin{array}{ll}
\sigma_1(\Xi,\Xi';h,f,x), &\mbox{ if $j=0$,}\\
\sigma_{2^j}(\Xi,\Xi';h,f,x)-\sigma_{2^{j-1}}(\Xi,\Xi';h,f,x), &\mbox{ if $j=1,2,\cdots$}
\end{array}\right.
\ee
We observe that since $h$ is a low pass filter, $\sigma_1$ is a filtered projection on the low frequency components for which $\lambda_k <1$. In the case when the $\lambda_k$'s are obtained as singular values of some operator, $\sigma_1$ clearly includes the projection on the kernel space of this operator.

In the case of analysis in $L^2(\mu^*)$, we may take $h$ to be the usual cut--off function: $h(t)=1$ if $0\le t<1$, $h(t)=0$ if $t\ge 1$. The operators $\sigma_n$ and $\tau_j$ are then the projections on $\Pi_n$ and the ``wavelet space'' $\Pi_{2^j}\ominus\Pi_{2^{j-1}}$ respectively. In the case of other $L^p$ spaces, for the purpose of analysis in uniform norm in particular, we need a smooth function $h$. A smooth function $h$ results in a well localized kernel for the operators.

Next, it is convenient to formulate some of our important assumptions in the form of a definition.

\begin{definition}\label{gaussupbdspace}
 The \textbf{heat kernel} on an admissible system  $\Xi=(\XX,d,\mu^*,\{\lambda_k\}, \{\phi_k\})$ is defined formally by
\be\label{heatkerndef}
K_t(\Xi; x,y)=\sum_{k=0}^\infty \exp(-\lambda_k^2t)\phi_k(x)\overline{\phi_k(y)}.
\ee 
We say that $\Xi$ satisfies the \textbf{Gaussian upper bound condition} (with exponent $q$) if there exists $q>0$, such that
\be\label{gaussupbdonxx}
|K_t(\Xi; x,y)|\le c_1t^{-q/2}\exp\left(-c_2\frac{d(x,y)^2}{t}\right), \qquad 0<t\le 1, \ x,y\in\XX.
\ee
In particular, the series in \eref{heatkerndef} converges for all $x,y\in\XX$.
\end{definition}

\begin{uda}\label{manifoldexample1}
{\rm (\textbf{the manifold case})
A typical example of a system satisfying the Gaussian upper bound condition is the following. Let $\XX$ be a compact ($C^\infty$), connected Riemannian manifold, $\mu^*$ be the Riemannian volume measure on $\XX$, and $d$ be the geodesic distance. We consider the Laplace--Beltrami operator on $\XX$, and let $\{-\lambda_k^2\}$ be the sequence of its eigenvalues, with corresponding sequence of eigenfunctions $\{\phi_k\}$. Then the system $(\XX, d, \mu^*, \{\lambda_k\}, \{\phi_k\})$ is an admissible system.
It is shown in \cite{kordyukov1991p} that under certain conditions on $\XX$, if $-\lambda_k^2$ are the eigenvalues of certain elliptic operators, and $\phi_k$ are the corresponding eigenfunctions, then \eref{gaussupbdonxx} is satisfied. Many other general situations in metric measure spaces are known \cite{grigor1997gaussian, grigor1999estimates, grigor2014heat} and references therein. \qed
}
\end{uda} 
\begin{rem}{\rm 
In the case when $\XX$ is a finite graph, the heat kernel corresponding to a vertex--based graph Laplacian satisfies the Gaussian upper bound condition for large values of $t$ rather than small values as in Definition~\ref{gaussupbdspace}. It is shown in \cite{friedman2004wave} that the solution of the wave equation corresponding to an edge--based Laplacian satisfies the property known as finite speed of wave propagation. In turn, this is equivalent to the Gaussian upper bound condition for small values of $t$ as stipulated in 
Definition~\ref{gaussupbdspace} \cite{sikora2004riesz, frankbern}.
}
\end{rem}

We will denote the ball of radius $r$ around $x$ by
\be\label{balldef}
\BB(x,r)=\{y\in\XX : d(x,y)\le r\},
\ee
and note that all balls are compact.

With this preparation, we are ready to state our wavelet--like representation theorem.

\begin{theorem}\label{digraphanaltheo}
Let $\Xi=(\XX, d, \mu^*, \{\lambda_k\}, \{\phi_k\})$ be an admissible system  satisfying the Gaussian upper bound condition with exponent $q>0$. Let $U$ be a partial isometry operator on $L^2(\mu^*)$, $\psi_k=U^*\phi_k$, and $\Xi'=(\XX, d, \mu^*, \{\lambda_k\}, \{\psi_k\})$ be the associated admissible system. We assume further that $\mu^*$ satisfies the regularity condition
\be\label{ballmeasurecond}
\mu^*(\BB(x,r))\le cr^q, \qquad x\in\XX, \ r>0.
\ee
Let $1\le p\le \infty$, $f\in L^1(\mu^*)+L^\infty(\mu^*)$,  $Uf\in X^p$, and $h$ be a sufficiently smooth low pass filter.\\
{\rm (a)} For $n\ge 1$, we have
\be\label{digraphgoodapprox}
E_{n/2,p}(\Xi;Uf)\le \|Uf-\sigma_n(\Xi,\Xi';h,f)\|_{\XX,\mu^*,p} \le cE_{n,p}(\Xi; Uf).
\ee
{\rm (b)} In the sense of $L^p(\mu^*)$ convergence,
\be\label{digraphlpseries}
Uf=\sum_{j=0}^\infty \tau_j(\Xi,\Xi';h,f).
\ee
{\rm (c)} If $p=2$, then
\be\label{digraphframebds}
\sum_{j=0}^\infty \|\tau_j(\Xi,\Xi';h,f)\|_{\XX,\mu^*,2}^2 \le \|Uf\|_{\XX,\mu^*,2}^2\le 5\sum_{j=0}^\infty \|\tau_j(\Xi,\Xi';h,f)\|_{\XX,\mu^*,2}^2.
\ee 
{\rm (d)} Let $0<\rho\le\infty$, $\gamma>0$. Then $Uf\in B_{p,\rho,\gamma}$ if and only if $\{\|\tau_j(\Xi,\Xi';h,f)\|_{\XX,\mu^*,p}\}_{j=0}^\infty\in \mathsf{b}_{\rho,\gamma}$.
\end{theorem}

\begin{rem}{\rm 
In the case when $U$ is the identity operator, Theorem~\ref{digraphanaltheo}
reduces to the corresponding theorems in \cite{mauropap} for the undirected case. \qed
}
\end{rem}

\bhag{Analysis on changing data}\label{twosystemsect}

Our constructions in this section are  motivated by the paper \cite{coifmanhirn} of Coifman and Hirn. In our terminology, they consider two admissible systems $\Xi_1=(\XX, d_1,\mu^*, \{\lambda_{1,j}\}, \{\phi_{1,j}\})$ and $\Xi_2=(\XX, d_2,\mu^*, \{\lambda_{2,k}\}, \{\phi_{2,k}\})$ with the base measure space $(\XX,\mu^*)$ being common to both, and $\{\phi_{1,j}\}$, $\{\phi_{2,k}\}$ being orthonormal systems. The diffusion distance between $x_1, x_2\in\XX$, with the points considered as elements of different copies of $\XX$ is defined by
$$
\|K_t(\Xi_1;x_1,\cdot)-K_t(\Xi_2;x_2,\cdot)\|_2=\left\{K_{2t}(\Xi_1; x_1,x_1) +K_{2t}(\Xi_2; x_2,x_2)-2\Re e\int_\XX K_t(\Xi_1;x_1,y)K_t(\Xi_2;y,x_2)d\mu^*(y)\right\}^{1/2}.
$$
Thus, the interaction term is the inner product expression above, which can be written in the form
$$
\Re e\sum_{j,k=0}^\infty \exp(-t(\lambda_{1,j}^2+\lambda_{2,k}^2))\overline{\langle \phi_{1,j}, \phi_{2,k}\rangle} \phi_{1,j}(x_1)\overline{\phi_{2,k}(x_2)},
$$
where $\langle \cdot,\cdot\rangle$ denotes in this context  the inner product in $L^2(\mu^*)$.

This prompts the following  more general set up. We are interested in the interaction of two systems $\Xi_1=(\XX_1,d_1,\mu_1^*,\{\lambda_{1,k}\}, \{\phi_{1,k}\})$
 and 
$\Xi_2=(\XX_2,d_2,\mu_2^*,\{\lambda_{2,k}\}, \{\phi_{2,k}\})$. 
We assume that a subset of ``landmarks'' $\YY_1\subset \XX_1$ is in a one--to--one correspondence with the corresponding subset of landmarks $\YY_2\subset \XX_2$, with some smoothness conditions on this correspondence as appropriate. To simplify notation, we will assume simply that $\YY_1=\YY_2=\YY\subseteq \XX_1\cap\XX_2$, and that there is a  positive measure $\nu^*$ supported on $\YY$ such that $\nu^*(\YY)=1$. The measure $\nu^*$ might well be discrete. 

Our first objective in this section is to examine which functions on $\XX_2$ can be lifted to functions on $\XX_1$ via the landmarks $\YY$, and how the smoothness of such a lifted function changes from that of the original function. One example is  the theory in Section~\ref{digraphharmanalsect}, 
where the ``lifted'' function is $Uf$, and the smoothness does not change.  A different model is the trace theorems for partial differential equations, where the smoothness of the trace/extension is typically different from the function one starts out with. 

As a simple starting point, we replace $\overline{\langle \phi_{1,j}, \phi_{2,k}\rangle}$ in the interaction term above by
\be\label{genconnectiondef}
\Gamma_{j,k}=\int_\YY \overline{\phi_{1,j}(y)}\phi_{2,k}(y)d\nu^*(y), \qquad j,k=0,1,\cdots.
\ee

Let $h$ be a smooth, low pass filter, and $\bs{\Gamma}$ be the matrix defined by \eref{genconnectiondef}. We define
\be\label{tensorsummkerneldef}
\Phi_{n,\otimes}(\Xi_1,\Xi_2;h, x_1,x_2)=\sum_{j, k=0}^\infty h\left(\frac{\lambda_{1,j}}{n}\right)h\left(\frac{\lambda_{2,k}}{n}\right)\Gamma_{j,k}\phi_{1,j}(x_1)\overline{\phi_{2,k}(x_2)}, \qquad x_1\in\XX_1, \ x_2\in\XX_2.
\ee
If $f\in L^1(\XX_2,\mu_2^*)+L^\infty(\XX_2,\mu_2^*)$, we define
\be\label{tensorsigmaopdef}
\sigma_{n,\otimes}(\Xi_1,\Xi_2;h,f,x_1)=\int_{\XX_2}f(x_2)\Phi_{n,\otimes}(\Xi_1,\Xi_2;h,x_1,x_2)d\mu_2^*(x_2), \qquad x_1\in\XX_1,\ n>0.
\ee 

In the case when $\Xi_1=(\XX, d_1,\mu^*, \{\lambda_{1,j}\}, \{\phi_{1,j}\})$ and $\Xi_2=(\XX, d_2,\mu^*, \{\lambda_{2,k}\}, \{\phi_{2,k}\})$, $\YY=\XX$, $\nu^*=\mu^*$, we have for each $k$,
$$
\sigma_{m,\otimes}(\Xi_1,\Xi_2;h,\phi_{2,k},x_1)=
h\left(\frac{\lambda_{2,k}}{m}\right)\sum_{j=0}^\infty h\left(\frac{\lambda_{1,j}}{m}\right)\Gamma_{j,k}\phi_{1,j}(x_1)
=h\left(\frac{\lambda_{2,k}}{m}\right)\sigma_m(\Xi_1,\Xi_1;h,\phi_{2,k}, x_1).
$$
It is not difficult to deduce from here,  under suitable conditions as in \cite{mauropap}, that 
 $\disp\lim_{m\to\infty}\sigma_{m,\otimes}(\Xi_1,\Xi_2;h,P)=P$ for every $P\in \Pi_\infty(\Xi_2)$. Moreover, the results in \cite{mauropap} show that if $P$ is sufficiently smooth in $\Xi_1$, then
$$
\|\sigma_{2^{m+1},\otimes}(\Xi_1,\Xi_2;h,P)-\sigma_{2^m,\otimes}(\Xi_1,\Xi_2;h,P)\|_{\XX_1,\mu_1^*,p}=\O(2^{-m\beta})
$$
for some $\beta>0$. In our more general case, this may or may not be true, depending upon the choice of $\YY$ and $\nu^*$. Nevertheless,
in order to obtain a theorem about the change in smoothness of a function when ``lifted'' from $\Xi_2$ to $\Xi_1$, we need to make a similar assumption (\eref{gammadecaycond} below).

\begin{theorem}\label{gentensortheo}
Let $\Xi_1=(\XX_1,d_1,\mu_1^*,\{\lambda_{1,k}\}, \{\phi_{1,k}\})$ and
 $\Xi_2=(\XX_2,d_2,\mu_2^*,\{\lambda_{2,k}\}, \{\phi_{2,k}\})$ be two admissible systems as in Definiton~\ref{systemdef}, each of which satisfies the Gaussian upper bound condition with exponents $q_1$, $q_2$ respectively, where $\{\phi_{1,j}\}$ and $\{\phi_{2,k}\}$ are orthonormalized with respect to $\mu_1^*$, respectively, $\mu_2^*$. Let $\YY$ be a measurable (with respect to both $\mu_1^*$ and $\mu_2^*$) subset of $\XX_1\cap \XX_2$, $\nu^*$ be a probability measure on $\YY$, $\Gamma_{j,k}$ be defined as in  \eref{genconnectiondef}. Let $h$ be a sufficiently smooth low pass filter, and $\mu_1^*$, $\mu_2^*$ satisfy the  following regularity condition  (cf. \eref{balldef}): 
 \be\label{ballmeasurecond2}
\mu_j^*(\BB_j(x,r))\le cr^{q_j}, \qquad x\in \XX_j, \ r>0, \ j=1,2.  
 \ee
Let $1\le p\le\infty$. We assume further that there exists a $\beta>0$ such that for any $m, n\ge 1$ and $P\in \Pi_n(\Xi_2)$,
\be\label{gammadecaycond}
 \|\sigma_{2^{m+1},\otimes}(\Xi_1,\Xi_2;h,P)-\sigma_{2^m,\otimes}(\Xi_1,\Xi_2;h,P)\|_{\XX_1,\mu_1^*,p}\le c2^{-m\beta}n^{q_1+(q_2-q_1)/p}\|P\|_{\XX_2,\mu_2^*,p}.
\ee
If $f\in L^p(\XX_2,\mu_2^*)$, and
\be\label{tenspseudobesovcond}
\sum_{m=0}^\infty 2^{m(q_1+(q_2-q_1)/p)}E_{2^m,p}(\Xi_2;f) <\infty,
\ee
then $\disp\mathcal{E}_\otimes(f)=\mathcal{E}_\otimes(\Xi_1,\Xi_2; f)=\lim_{n\to\infty}\sigma_{2^n,\otimes}(\Xi_1,\Xi_2;h,f)$ exists in the sense of $L^p(\Xi_1,\mu_1^*)$, and
\be\label{tensgoodapprox}
E_{2^n,p}(\Xi_1;\mathcal{E}_\otimes(f))\le c\left\{\sum_{m=n}^\infty 2^{m(q_1+(q_2-q_1)/p)}E_{2^m,p}(\Xi_2;f)+2^{n(q_1+(q_2-q_1)/p-\beta)}\|f\|_{\Xi_2,p}\right\}.
\ee
In particular, if $\beta>\gamma>q_1+(q_2-q_1)/p$, $0<\rho\le\infty$, and $f\in B_{p,\rho,\gamma}(\Xi_2)$, then $\mathcal{E}_\otimes(f)\in B_{p,\rho,\gamma-q_1-(q_2-q_1)/p}$.
\end{theorem}

We note that even if $\Xi_1=\Xi_2$ in the above theorem, a straighforward application of the theorem leads to a  loss of smoothness in $f$. This is a technicality due in part to the very general nature of $\YY$ and $\nu^*$, but also due to the definition of $\sigma_{n,\otimes}$. It is possible to eliminate  this anomaly by defining the analogue of the various quantities in one system theory in our context in an abstract manner rather than assuming the identification of some elements of $\XX_1$ with those in $\XX_2$ via $\YY$.

In particular, we want to define a version of the Gaussian upper bound condition appropriate to two systems considered jointly. To motivate this discussion, we first make an observation, continuing the same scenario as in Theorem~\ref{gentensortheo}. 
\begin{uda}\label{jointgaussexample}
{\rm
 Let $\Xi_1=(\XX_1,d_1,\mu_1^*,\{\lambda_{1,j}\}, \{\phi_{1,j}\})$ and
 $\Xi_2=(\XX_2,d_2,\mu_2^*,\{\lambda_{2,k}\}, \{\phi_{2,k}\})$ be two admissible systems as in Definiton~\ref{systemdef}, each of which satisfies the Gaussian upper bound condition with exponents $q_1$, $q_2$ respectively, where $\{\phi_{1,j}\}$ and $\{\phi_{2,k}\}$ are orthonormalized with respect to $\mu_1^*$, respectively, $\mu_2^*$. Let $\YY$ be a measurable (with respect to both $\mu_1^*$ and $\mu_2^*$) subset of $\XX_1\cap \XX_2$, $\nu^*$ be a probability measure on $\YY$. We define, in this example only,
\be\label{specialjointdistancedef}
d_{1,2}(x_1,x_2)=\inf_{y\in\YY}\left(d_1(x_1,y)+d_2(y,x_2)\right).
\ee
Since $d_{1,2}$ is defined on $\XX_1\times\XX_2$, it cannot be called a distance in the usual sense of the term, but it is clearly non--negative, satisfies the symmetry condition $d_{1,2}(x_1,x_2)=d_{2,1}(x_2,x_1)$, and the analogues of the triangle inequality:
\bea\label{specialjointtriangineq}
d_{1,2}(x_1,x_2)&\le& c\left(d_1(x_1, x_1')+d_{1,2}(x_1',x_2)\right), \qquad x_1,x_1'\in\XX_1, \ x_2\in\XX_2,\nonumber\\
d_{1,2}(x_1,x_2)&\le& c\left(d_{1,2}(x_1,x_2')+d_2(x_2', x_2)\right), \qquad x_1\in\XX_1, \ x_2,x_2'\in\XX_2.
\eea
Next, we let $\ell_{j,k}=(\lambda_{1,j}^2+\lambda_{2,k}^2)^{1/2}$, and
\be\label{specialgenjointheatkernel}
K_t(\Xi_1,\Xi_2;x_1,x_2)=\sum_{k, j=0}^\infty  \exp(-\ell_{j,k}^2t)\Gamma_{j,k}\phi_{1,j}(x_1)\overline{\phi_{2,k}(x_2)}.
\ee
An application of \eref{gaussupbdonxx} shows that
\bea\label{genjointheatkernest}
|K_t(\Xi_1,\Xi_2;x_1,x_2)|
&=& \left|\int_\YY K_t(\Xi_1; x_1,y)\overline{K_t(\Xi_2; x_2,y)}d\nu^*(y)\right|\nonumber\\
&\le&\int_\YY |K_t(\Xi_1; x_1,y)||K_t(\Xi_2; x_2,y)|d\nu^*(y)\nonumber\\
&\le& ct^{-(q_1+q_2)/2}\int_\YY \exp\left(-c\frac{d_1(x_1,y)^2+d_2(y,x_2)^2}{t}\right)d\nu^*(y)\nonumber\\
&\le& ct^{-(q_1+q_2)/2}\exp\left(-c\frac{d_{1,2}(x_1,x_2)^2}{t}\right).
\eea
\qed
}
\end{uda}

The next example is meant to illustrate the need for a more general set up which arises in a very routine extension of functions on the closed unit disc of $\RR^2$ to the Euclidean sphere $\SS^2$ in $\RR^3$.

\begin{uda}\label{jacobiexample}
{\rm
Let $\XX_1$ be the Euclidean hemisphere
$$
\XX_1=\{\textbf{p}(\theta,\phi)=(\sin\theta\cos\phi,\sin\theta\sin\phi, \cos\theta) : \theta\in [0,\pi/2], \ \phi\in (-\pi,\pi]\},
$$
$d_1$ be the geodesic distance on $\XX_1$, $\mu_1^*$ be the area measure, normalized to be a probability measure, $\lambda_{1,(j,\ell)}=\sqrt{(j+\ell)(j+\ell+1)}$, $j=0,1,\cdots$, $\ell\in\ZZ$, and
$$
\phi_{1,(j,\ell)}(\textbf{p}(\theta,\phi))= \sqrt{2}\sin^{|\ell|}\theta p_{2j+1}^{(|\ell|,|\ell|)}(\cos\theta)\exp(i\ell\phi),
$$
where, in this example, $p_k^{(\a,\beta)}$ denotes the orthonormalized Jacobi polynomial of degree $k$; i.e., a polynomial of degree $k$ with positive leading coefficients, such that for integer $k, j\ge 0$,
$$
\int_{-1}^1  p_k^{(\a,\beta)}(x)p_m^{(\a,\beta)}(x)(1-x)^\a(1+x)^\beta dx=
\left\{\begin{array}{ll}
1, & \mbox{if $k=m$},\\
0, &\mbox{otherwise}.
\end{array}\right.
$$
For $\XX_2$, we take the unit disc $\{(x,y)\in\RR^2: x^2+y^2\le1\}$, $d_2$ to be the Euclidean distance, $\mu_2^*$ to be the area measure, let $\textbf{q}(\theta,\phi)=(\sin\theta\cos\phi, \sin\theta\sin\phi)$, $\theta\in [0,\pi/2]$, $\phi\in (-\pi,\pi]$, $\lambda_{2,(k,m)}=k|m|+k+|m|$,  and define $\phi_{2,(k,m)}$ by
$$
\phi_{2,(k,m)}(\textbf{q}(\theta,\phi))=\sqrt{\frac{2^{|m|+2}}{\pi}}\sin^{|m|}\theta\cos\theta p_k^{(|m|,1)}(\cos(2\theta))\exp(im\phi), \qquad m\in\ZZ,\ k=0,1,\cdots.
$$
In this example, each member of the systems $\{\phi_{1,(j,\ell)}\}$, $\{\phi_{2,(k,m)}\}$ is equal to $0$ on $\XX_1\cap\XX_2$. So, the construction as in Example~\ref{jointgaussexample} does not work. There is, of course, an obvious one--to--one correspondence between each point $\textbf{p}(\theta,\phi)\in\XX_1$ and $\textbf{q}(\theta,\phi)\in\XX_2$, but since the measures are different, we are in a different setting from the one described in \cite{coifmanhirn} as well. 
As an interaction term, we find it convenient to define 
for $j, k=0,1,\cdots$, $\ell, m\in\ZZ$,
$$
A_{(j,\ell), (k,m)}=\left\{\begin{array}{ll}
\disp\sqrt{\pi}2^{-|m|-5/4}a_{j,k}^{(m)}, &\mbox{if $\ell=m$, $0\le k\le j$},\\
0, &\mbox{otherwise},
\end{array}\right.
$$
where
$$
a_{j,k}^{(m)}=\int_{-1}^1p_j^{(|m|,1/2)}(x)p_k^{(|m|,1)}(x)(1-x)^{|m|}(1+x)dx, \qquad 0\le k\le j,\ j=0,1, \cdots, \ m\in \ZZ.
$$
Explicit expressions for $a_{j,k}^{(m)}$ are known, cf. \cite[Lecture~7]{askeyspecial}. Using the fact  (\cite[Theorem~4.1, Formula~(4.3.4)]{szego}) that for $\a>-1$, $\theta\in [0,\pi/2]$,
\be\label{jacobiultra}
p_{2j+1}^{(\a,\a)}(\cos\theta)=2^{\a/2+3/4}\cos\theta p_j^{(\a,1/2)}(\cos(2\theta)), \qquad j=0,1,\cdots, 
\ee
we now deduce that
$$
\sum_{k=0}^\infty\sum_{m\in\ZZ}A_{(j,\ell),(k,m)}\phi_{2,(k,m)}(\textbf{q}(\theta_2,\phi_2))= \phi_{1,(j,\ell)}((\textbf{p}(\theta_2,\phi_2)).
$$
Since $\phi_{1,(j,\ell)}$ are Dirichlet eigenfunctions of the Laplace--Beltrami operator on $\XX_1$ corresponding to the eigenvalue $-\lambda_{1,(j,\ell)}^2$, it follows from the results in \cite{kordyukov1991p, davies1990heat} that
\begin{eqnarray*}
\lefteqn{\left|\sum_{j, k=0}^\infty\sum_{m,\ell\in\ZZ}\exp(-(j+\ell)(j+\ell+1)t)A_{(j,\ell),(k,m)}\phi_{1,(j,\ell)}((\textbf{p}(\theta,\phi))\overline{\phi_{2,(k,m)}(\textbf{q}(\theta_2,\phi_2))}\right|}\\
&=& \left|\sum_{j=0}^\infty\sum_{\ell\in\ZZ}\exp(-(j+\ell)(j+\ell+1)t)\phi_{1,(j,\ell)}((\textbf{p}(\theta_1,\phi_1))\overline{\phi_{1,(j,\ell)}(\textbf{p}(\theta_2,\phi_2)})\right|\\
&\le& c_1t^{-1}\exp\left(-c_2\frac{d_1((\textbf{p}(\theta_1,\phi_1), (\textbf{p}(\theta_2,\phi_2))^2}{t}\right).
\end{eqnarray*}
We note finally that setting 
$$
d_{1,2}(\textbf{p}(\theta_1,\phi_1), \textbf{q}(\theta_2,\phi_2))=d_1((\textbf{p}(\theta_1,\phi_1), (\textbf{p}(\theta_2,\phi_2)),
$$
the function $d_{1,2} : \XX_1\times\XX_2$ satisfies a triangle inequality analogous to \eref{specialjointtriangineq}.
\qed
}
\end{uda}
Motivated by these examples, we make first the next definition of a joint distance.
\begin{definition}\label{jointdistdef}
Let $(\XX_1,d_1)$, $(\XX_2,d_2)$ be quasi--metric spaces. A function $d_{1,2} :\XX_1\times\XX_2\to [0,\infty)$, will be called a \textbf{joint distance} for the pair of quasi--metric spaces if
the following triangle inequalities are satisfied. 
\bea\label{jointtriangineq}
d_{1,2}(x_1,x_2)&\le& c\left(d_1(x_1, x_1')+d_{1,2}(x_1',x_2)\right), \qquad x_1,x_1'\in\XX_1, \ x_2\in\XX_2,\nonumber\\
d_{1,2}(x_1,x_2)&\le& c\left(d_{1,2}(x_1,x_2')+d_2(x_2', x_2)\right), \qquad x_1\in\XX_1, \ x_2,x_2'\in\XX_2.
\eea
We define
\be\label{jointsymmetry}
d_{2,1}(x_2,x_1)=d_{1,2}(x_1,x_2), \qquad x_1\in \XX_1,\  x_2\in\XX_2,
\ee 
\end{definition}
For $r>0$, $x_1\in\XX_1$, $x_2\in\XX_2$, we define the balls
\bea\label{jointballdef}
\BB_1(x_1,r)=\{z\in \XX_1 : d_1(x_1,z)\le r\}, &\quad& \BB_2(x_2,r)=\{z\in \XX_2 : d_2( x_2,z)\le r\}.\nonumber\\
\BB_{1,2}(x_1,r)=\{z\in \XX_2 : d_{1,2}(x_1,z)\le r\}, &\quad& \BB_{2,1}(x_2,r)=\{z\in \XX_1 : d_{2,1}( x_2,z)\le r\}.
\eea

Next, we would like to allow more flexible connections than that provided by the matrix $\Gamma$ as in \eref{genconnectiondef}.    Example~\ref{jacobiexample} demonstrates the need to do so. Another motivation is to incorporate some probablistic transition from a point $x_1\in\XX_1$ to $x_2\in\XX_2$. These considerations prompt  the following very general definition analogous to Definition~\ref{gaussupbdspace}.

\begin{definition}\label{jointgaussbd}
Let $\Xi_1=(\XX_1,d_1,\mu_1^*,\{\lambda_{1,k}\}, \{\phi_{1,k}\})$ and $\Xi_2=(\XX_2,d_2,\mu_2^*,\{\lambda_{2,k}\}, \{\phi_{2,k}\})$ be two admissible systems as in Definiton~\ref{systemdef}, $d_{1,2} $ be the joint distance between the two systems. Let $\mathbf{A}=(A_{j,k})_{j,k=0}^\infty$,
$\bs{\lambda}=(\ell_{j,k})_{j,k=0}^\infty$,
where each $\ell_{j,k}\ge 0$, and for each $u>0$, the set $\{(j,k) : \ell_{j,k}<u\}$ is finite. We define formally for $x_1\in\XX_1$, $x_2\in\XX_2$, $t>0$, the \textbf{joint heat kernel} by
\bea\label{jointheatkerndef}
K_t(\Xi_1,\Xi_2;x_1,x_2)&=&K_t(\Xi_1,\Xi_2;\mathbf{A}, \bs{\lambda}; x_1,x_2)=\sum_{j,k=0}^\infty \exp(-\ell_{j,k}^2 t)A_{j,k}\phi_{1,j}(x_1)\overline{\phi_{2,k}(x_2)}\nonumber\\
&=&\lim_{n\to\infty}\sum_{j,k:\ell_{j,k}<n} \exp(-\ell_{j,k}^2 t)A_{j,k}\phi_{1,j}(x_1)\overline{\phi_{2,k}(x_2)}.
\eea
We say that $\Xi_1$, $\Xi_2$ satisfy The \textbf{joint Gaussian upper bound condition} (with exponent $Q>0$) if 
\be\label{jointondiagbd}
\sum_{j,k : \ell_{j,k}<n}|A_{j,k}\phi_{1,j}(x_1)\overline{\phi_{2,k}(x_2)}|\le cn^Q, \ n\ge 2,
\ee
the limit in \eref{jointheatkerndef} exists for all $x_1\in\XX_1$, $x_2\in\XX_2$, and
\be\label{jointoffdiagbd}
|K_t(\Xi_1,\Xi_2;x_1,x_2)| \le c_1t^{-c}\exp\left(-c_2\frac{d_{1,2}(x_1,x_2)^2}{t}\right), 
\qquad x_1\in\XX_1,\ x_2\in\XX_2.
\ee
We will refer to $\mathbf{A}$ as the \textbf{connection matrix}, its elements as \textbf{connection coefficients}, and the numbers $\ell_{j,k}$ as \textbf{joint eigenvalues}.
\end{definition}

\begin{rem}{
\rm
When $\Xi_1=\Xi_2$, $\mathbf{A}$ is the identity matrix, $\ell_{j,k}=\lambda_k\delta_{j,k}$, then \eref{jointondiagbd} takes the form
$$
\sum_{k : \lambda_k <n}\exp(-\lambda_k^2t)|\phi_{1,k}(x_1)|^2\le cn^Q, \qquad x_1\in \XX_1,\ n\ge 2.
$$
We have proved in \cite[Proposition~4.1]{frankbern} that this is equivalent to the statement that
$$
K_t(\Xi_1,\Xi_1;x_1,x_1)\le ct^{-Q/2}, \qquad 0<t<1.
$$
Thus, the joint Gaussian upper bound condition in Definition~\ref{jointgaussbd} reduces to the Gaussian upper bound condition for a single system defined in Definition~\ref{gaussupbdspace}. \qed
}
\end{rem}
\begin{uda}\label{bernsteinexample}
{\rm
In the manifold case as in Example~\ref{manifoldexample1}, if  $\psi_k$ is a derivative of $\phi_k$, then an integration by parts argument shows that $\{\psi_k\}$ is a Bessel system. It is known that under certain conditions, the systems $\Xi_1=(\XX,d,\mu^*\!,\{\lambda_k\}, \{\phi_k\})$, $\Xi_2=(\XX,d,\mu^*,\{\lambda_k\}, \{\psi_k\})$ satisfy the joint Gaussian upper bound condition (\cite{grigoryan1995upper} and references therein). \qed}
\end{uda}

In the remainder of this section, let $\Xi_1=(\XX_1,d_1,\mu_1^*,\{\lambda_{1,k}\}, \{\phi_{1,k}\})$ and $\Xi_2=(\XX_2,d_2,\mu_2^*,\{\lambda_{2,k}\}, \{\phi_{2,k}\})$ be two admissible systems satisfying the joint Gaussian upper bound condition with exponent $Q>0$, and $d_{1,2}$, $\mathbf{A}$, $\bs{\lambda}$ be   as in Definition~\ref{jointgaussbd}.

Next, we define the analogues of the operators $\sigma_n$.

If $h$ is a low pass filter, we define the kernel $\Phi_n$ for $x_1\in\XX_1, \ x_2\in\XX_2$ by
\be\label{jointsummkerneldef}
\Phi_n(\Xi_1,\Xi_2; h, x_1,x_2)=\Phi_n(\Xi_1,\Xi_2;\mathbf{A}, \bs{\lambda}; h, x_1,x_2)=\sum_{j, k=0}^\infty h\left(\frac{\ell_{j,k}}{n}\right)A_{j,k}\phi_{1,j}(x_1)\overline{\phi_{2,k}(x_2)}.
\ee
If $f\in L^1(\XX_2,\mu_2^*)+L^\infty(\XX_2,\mu_2^*)$, $x_1\in\XX_1$, $n>0$, we define
\bea\label{jointsigmaopdef}
\sigma_n(\Xi_1,\Xi_2; h,f,x_1)&=&\sigma_n(\Xi_1,\Xi_2;\mathbf{A}, \bs{\lambda}; h,f,x_1)=
\int_{\XX_2}f(x_2)\Phi_n(\Xi_1,\Xi_2;h,x_1,x_2)d\mu_2^*(x_2)\nonumber\\&=& \sum_{j,k=0}^\infty h\left(\frac{\ell_{j,k}}{n}\right)A_{j,k}\hat{f}(\Xi_2;k)\phi_{1,j}(x_1).
\eea
\begin{rem}\label{diffhilbertrem}{\rm
In the case when the matrix $\mathbf{A}$ is the identity matrix, $\{\phi_{1,j}\}$, $\{\phi_{2,k}\}$ are orthonormal systems, $U$ is a unitary operator from $L^2(\XX_2,\mu_2^*)$ into $L^2(\XX_1,\mu_1^*)$  such that $\phi_{2,k}=U^*\phi_{1,k}$ for each $k$, then the kernels $\Phi_n$ and the corresponding operators $\sigma_n$ correspond to these quantities defined in 
Section~\ref{digraphharmanalsect} in the case when $\XX_1=\XX_2$, $\mu_1^*=\mu_2^*$. Because of the special relationship between the two systems, the results in Section~\ref{digraphharmanalsect} are valid in this more general context, assuming the Gaussian upper bound condition only for the system $\{\phi_{1,j}\}$, and not the system  $\{\phi_{2,k}\}$. In the greater generality of this section, where there is no such natural connection between the systems, we need to assume a joint Gaussian upper bound condition for the following theorems to hold. \qed
}
\end{rem}

Our generalization of Theorem~\ref{gentensortheo} is the following.

\begin{theorem}\label{extensiontheo}
Let $\Xi_1=(\XX_1,d_1,\mu_1^*,\{\lambda_{1,k}\}, \{\phi_{1,k}\})$ and $\Xi_2=(\XX_2,d_2,\mu_2^*,\{\lambda_{2,k}\}, \{\phi_{2,k}\})$ be two admissible systems satisfying the joint Gaussian upper bound condition with exponent $Q>0$, and $\mathbf{A}$, $\bs{\lambda}$   as in Definition~\ref{jointgaussbd}. Let $h$ be a sufficiently smooth low pass filter, and $\mu_1^*, \mu_2^*$ satisfy the \textbf{ joint regularity condition}
\be\label{ballmeasurecond1}
\mu_1^*(\BB_{2,1}(x_2,r))\le cr^{q_1}, \quad \mu_2^*(\BB_{1,2}(x_1,r))\le cr^{q_2}, \qquad x_1\in\XX_1,\ x_2\in\XX_2,\  r>0,
\ee
for some constants $c, q_1, q_2>0$ independent of $x_1$ and $r$. Let $1\le p\le \infty$, and $f\in L^p(\XX_2,\mu_2^*)$. We assume further that there exists $\beta>Q-q_2-(q_1-q_2)/p$ such that for any $n\ge 1$ and $P\in \Pi_n(\Xi_2)$,
\be\label{adecaycond}
\|\sigma_{2^{m+1}}(\Xi_1,\Xi_2;\mathbf{A}, \bs{\lambda};h,P)-\sigma_{2^m}(\Xi_1,\Xi_2;\mathbf{A}, \bs{\lambda};h,P)\|_{\XX_1,\mu_1^*,p}\le c2^{-m\beta} n^{Q-q_2-(q_1-q_2)/p}\|P\|_{\XX_2,\mu_2^*,p}, \qquad m=1,2,\cdots.
\ee
If 
\be\label{pseudobesovcond}
\sum_{m=0}^\infty 2^{m(Q-q_2-(q_1-q_2)/p)}E_{2^m,p}(\Xi_2;f)<\infty,
\ee
 then $\lim_{n\to\infty} \sigma_{2^n}(\Xi_1,\Xi_2;\mathbf{A}, \bs{\lambda}; h,f)=\mathcal{E}(f)=\mathcal{E}(\Xi_1,\Xi_2;\mathbf{A}, \bs{\lambda};f)$ exists in the sense of $L^p(\XX_1,\mu_1^*)$ and
\bea\label{twospacegoodapprox}
\lefteqn{\|\mathcal{E}(f)-\sigma_{2^n}(\Xi_1,\Xi_2;\mathbf{A}, \bs{\lambda};h,f)\|_{\XX_1,\mu^*,p}}\nonumber\\
&\le& c\left\{\sum_{m=n}^\infty 2^{m(Q-q_2-(q_1-q_2)/p)}E_{2^m,p}(\Xi_2;f) + 2^{n(Q-q_2-(q_1-q_2)/p-\beta)}\|f\|_{\XX_2,\mu_2^*,p}\right\}.
\eea
In particular, if $\a>0$, and  $\a\ell_{j,k}\ge \lambda_{1,j}$ for all $j, k=0,1,\cdots$, then $\sigma_n(\Xi_1,\Xi_2;\mathbf{A}, \bs{\lambda};h,f)\in \Pi_{\a n}(\Xi_1)$. Further, if   $Q-q_2-(q_1-q_2)/p<\gamma<\beta$, and $f\in B_{p,\rho, \gamma}(\Xi_2)$ then $\mathcal{E}(f)\in B_{p,\rho, \gamma-Q+q_2+(q_1-q_2)/p}(\Xi_1)$.
\end{theorem}

\begin{rem}{\rm
In the case when $p=2$, the condition \eref{adecaycond} can be replaced by the condition that the spectral norms of dyadic sectors $(A_{j,k})_{2^{m-1}\le \ell_{j,k}<2^m,\ \lambda_{2,k}<n}$
are bounded from above by $c2^{-m\beta}n^{Q-q_2-(q_1-q_2)/p}$. 
\qed
}
\end{rem}
\begin{rem}{\rm
In Theorem~\ref{extensiontheo}, we do not assume that $Q-q_2-(q_1-q_2)/p>0$. Therefore, depending upon the sign of this expression, $\mathcal{E}$ is a smoothing operator analogous to an integral operator, or smoothness preserving, or coarsening operation analogous to a differential operator. \qed
}
\end{rem}
\begin{rem}{\rm
The exponent in place of $Q-q_2-(q_1-q_2)/p$   in the ``reverse'' direction is $Q-q_1-(q_2-q_1)/p$. If $Q\not=(q_1+q_2)/2$, this is not equal to the negative of the first exponent. 
 This leads
to the fact that it may not hold in general that under the conditions of 
Theorem~\ref{extensiontheo}, $f\in B_{p,\rho, \gamma}(\Xi_2)$ \textbf{if and only if } $\mathcal{E}(f)\in B_{p,\rho, \gamma-Q+q_2+(q_2-q_1)/p}(\Xi_1)$. Moreover, in order to obtain $f$ back from $\mathcal{E}(f)$, one needs a construction with  $\mathbf{A}$ replaced by its ``inverse'' in some sense, and with some restrictions on $\bs{\lambda}$. We do not find it worthwhile to go into these details;  this would lead only to a more complicated re--statement of Theorem~\ref{extensiontheo} in two directions, without adding any new concepts.  One sided analogues of Theorem~\ref{digraphanaltheo} are quite easy to prove though, and we omit these details as well. \qed
}
\end{rem}

\begin{rem}{\rm
In the scenario considered in Example~\ref{jacobiexample}, $Q=q_1=q_2=2$, and $\beta$ can be chosen arbitrarily large. If $f: \XX_2\to\CC$, The operator $\mathcal{E}(f)$ is the usual lifting to the sphere:
$$
\mathcal{E}(f)(\textbf{p}(\theta,\phi))=f(\textbf{q}(\theta,\phi)),
$$
without requiring \eref{pseudobesovcond}. Theorem~\ref{extensiontheo} then confirms the well known fact that $\mathcal{E}$ preserves the smoothness of the function $f$.
\qed
}
\end{rem}

\bhag{Proofs}\label{pfsect}

\subsection{Proofs of the theorems in Section~\ref{digraphharmanalsect} and Section~\ref{twosystemsect}.}\label{decosect}

\noindent
\textsc{Proof of Theorem~\ref{digraphanaltheo}.}
As we mentioned earlier,  this theorem is a direct consequence of the fact that
$$
\hat{f}(\Xi';k)=\int_\XX f(y)\overline{\psi_k(y)}d\mu^*(y)=\int_\XX f(y)\overline{(U^*\phi_k)(y)}d\mu^*(y)=\widehat{Uf}(\Xi;k),
$$
so that
$$
\sigma_n(\Xi,\Xi';h,f)=\sigma_n(\Xi,\Xi;h,Uf).
$$
The statements (a), (b), (c) of this theorem follow from \cite[Theorem~2.1]{mauropap}, and the statement (d) follows from \cite[Theorem~3.1]{mauropap}, applying these theorems with $Uf$ in place of $f$. 
In \cite{mauropap}, the conditions on the system $\Xi$ were stated in terms of upper bounds on $\sum_{k: \lambda_k<n}|\phi_k(x)|^2$ and the finite speed of wave propagation. 
These were shown in \cite{frankbern} to be equivalent to the Gaussian upper bound condition. (See also Theorem~\ref{sikoratheo} below.) \qed

It is convenient to prove Theorem~\ref{extensiontheo} first, before Theorem~\ref{gentensortheo}. An important ingredient of our proof is the following theorem.

\begin{theorem}\label{mainappltheo}
Let $\Xi_1=(\XX_1,d_1,\mu_1^*,\{\lambda_{1,k}\}, \{\phi_{1,k}\})$ and $\Xi_2=(\XX_2,d_2,\mu_2^*,\{\lambda_{2,k}\}, \{\phi_{2,k}\})$ be two admissible systems satisfying the joint Gaussian upper bound condition with exponent $Q>0$, and $\mathbf{A}$, $\bs{\lambda}$   as in Definition~\ref{jointgaussbd}. Let $h$ be a sufficiently smooth low pass filter, and $\mu_1^*, \mu_2^*$ satisfy the  joint regularity condition \eref{ballmeasurecond1}.
 Let $1\le p\le \infty$, and $f\in L^p(\XX_2,\mu_2^*)$. Then for $n\ge 1$,
\be\label{uniformsigmabd}
\|\sigma_n(\Xi_1,\Xi_2;h,f)\|_{\XX_1,\mu_1^*,p}\le cn^{Q-q_2-(q_1-q_2)/p}\|f\|_{\XX_2,\mu_2^*,p}, \qquad f\in L^p(\XX_2,\mu_2^*).
\ee
\end{theorem}

\begin{uda}\label{mauropaprem}
{\rm Taking $\Xi_1=\Xi_2$, $Q=q_1=q_2$, $\mathbf{A}$ to be the identity matrix, $\bs{\lambda}=\{\lambda_j\delta_{j,k}\}$, and $\{\phi_{1,j}\}$ to be an orthonormal system, we recover \cite[Proposition~2.1]{mauropap}. This, in turn, leads to the other theorems in \cite{mauropap} used in the proof of Theorem~\ref{digraphanaltheo}. 
}\qed
\end{uda}

\begin{uda}\label{frankrem}
{\rm We consider an admissible system $\Xi_1=(\XX, d,\mu^*, \{\lambda_k\}, \{\phi_k\})$ as in Example~\ref{bernsteinexample} for $\XX_2$, define $\psi_k$ as in that example, and let $\Xi_2=(\XX, d,\mu^*, \{\lambda_k\}, \{\psi_k\})$. Then Theorem~\ref{mainappltheo} leads, in particular, to \cite[Theorem~2.1]{frankbern}. 
}\qed 
\end{uda}

We will prove Theorem~\ref{mainappltheo} by using the Tauberian theorem, Theorem~\ref{maintaubertheo}, to be proved in Sub--section~\ref{taubsect}. For reasons of organizational convenience, we assume this theorem in this section, and note that if $h$ is a sufficiently smooth low pass filter, then it can be extended to $\RR$ as an even function, that satisfies the conditions on $H$ stipulated in Theorem~\ref{maintaubertheo}.

In order to prove Theorem~\ref{mainappltheo}, we first prove a lemma.

\begin{lemma}\label{fundalemma}
 If $\mu_2^*$ satisfies  \eref{ballmeasurecond1}, then for $S>q_2$,
\be\label{fundaineq}
\int_{\XX_2}\frac{d\mu_2^*(x_2)}{\max(1, (nd_{1,2}(x_1, x_2)))^S} \le cn^{-q_2}.
\ee
\end{lemma}

\begin{Proof}\ 
In this proof only, let $\mathcal{A}_0=\BB_{1,2}(x_1,1/n)$ and $\mathcal{A}_m=\BB_{1,2}(x_1, 2^{m}/n)\setminus \BB_{1,2}(x_1, 2^{m-1}/n)$, $m=1,2,\cdots$. Then \eref{ballmeasurecond1} shows that
$$
\mu_2^*(\mathcal{A}_m)\le c(2^m/n)^{q_2}, \qquad m=0,1,\cdots.
$$
Consequently,
$$
\int_{\XX_2}\frac{d\mu_2^*(x_2)}{\max(1, (nd_{1,2}(x_1, x_2)))^S}= \sum_{m=0}^\infty \int_{\mathcal{A}_m}\frac{d\mu_2^*(x_2)}{\max(1, (nd_{1,2}(x_1, x_2)))^S}
\le cn^{-q_2}\sum_{m=0}^\infty 2^{-m(S-q_2)}\le cn^{-q_2}.
$$
\end{Proof}

In the remainder of this section, we assume that $h$ is sufficiently smooth so as to satisfy the conditions \eref{Hbvcondnew}--\eref{Hbvintbdnew}
in Theorem~\ref{maintaubertheo} with  $S>\max(q_1,q_2, Q+1)$.

\noindent
\textsc{Proof of Theorem~\ref{mainappltheo}.} In this proof,  we denote $\Phi_n(\Xi_1,\Xi_2;h,x_1,x_2)$ by $\Phi_n(x_1,x_2)$.
For each $x_1\in\XX_1$ and $x_2\in\XX_2$, we define the measure $\mu$ on $[0,\infty)$ by
\be\label{pf9eqn1}
\mu([0,u))=\mu(x_1,x_2;[0, u))=\sum_{j,k :\ell_{j,k}<u}A_{j,k}\phi_{1,j}(x_1)\overline{\phi_{2,k}(x_2)}, \qquad u\ge 0.
\ee
Using \eref{jointondiagbd} and the fact that each of the functions $\phi_{1,j}$, $\phi_{2,k}$ are bounded on their domains, it follows that (cf. \eref{muchristbd})
\be\label{pf9eqn2}
\tn\mu\tn_Q \le c, \qquad x_1\in\XX_1,\ x_2\in\XX_2.
\ee
The estimate \eref{jointoffdiagbd} is the same as \eref{muheatgaussbd} with $d_{1,2}(x_1,x_2)$ in place of $r$. Hence, the conditions of Theorem~\ref{maintaubertheo} are satisfied. 
Therefore, Theorem~\ref{maintaubertheo} shows that for each $x_1\in\XX_1$, $x_2\in\XX_2$,
\be\label{pf9eqn3}
\left|\int_0^\infty h(u/n)d\mu(x_1,x_2;\mu)\right|=\left|\sum_{j, k=0}^\infty h(\ell_{j,k}/n)A_{j,k}\phi_{1,j}(x_1)\overline{\phi_{2,k}(x_2)}\right| =|\Phi_n(x_1,x_2)|\le c\frac{n^Q}{\max(1, (nd_{1,2}(x_1, x_2)))^S}.
\ee
Therefore, Lemma~\ref{fundalemma} (applied once with $\mu_1^*$ and once with $\mu_2^*$ as stated) leads to
$$
\int_{\XX_1}|\Phi_n(x_1,x_2)|d\mu_1^*(x_1)\le cn^{Q-q_1}, \qquad x_2\in\XX_2,
$$
and
$$ \int_{\XX_2}|\Phi_n(x_1,x_2)|d\mu_2^*(x_2)\le cn^{Q-q_2}, \qquad x_1\in\XX_1.
$$
These inequalities lead to \eref{uniformsigmabd} for $p=1$ and $p=\infty$ respectively. The general case follows by an application of the Riesz--Thorin interpolation theorem. 
\qed

\noindent
\textsc{Proof of Theorem~\ref{extensiontheo}.} In this proof, we write $\sigma_n(f)$ in place of $\sigma_n(\Xi_1,\Xi_2;h,f)$, and $q$ in place of $Q-q_2-(q_1-q_2)/p$. Let $P\in \Pi_n(\Xi_2)$. Then \eref{adecaycond}  implies that
\be\label{pf3eqn1}
\sum_{m=0}^\infty \|\sigma_{2^{m+1}}(P)-\sigma_{2^m}(P)\|_{\Xi_1,\mu_1^*,p}\le c\sum_{m=0}^\infty 2^{-m\beta}n^q\|P\|_{\Xi_2,\mu_2^*,p}\le cn^q\|P\|_{\Xi_2,\mu_2^*,p}.
\ee
Therefore, there exists $\mathcal{E}(P)\in L^p(\XX_1, \mu_1^*)$ such that
\be\label{pf3eqn2}
\mathcal{E}(P)=\sigma_1(P)+\sum_{m=0}^\infty \left(\sigma_{2^{m+1}}(P)-\sigma_{2^m}(P)\right),
\ee
where the series converges in the sense of $L^p(\XX_1,\mu_1^*)$. 
Moreover, $\mathcal{E}$ is a linear operation on $\Pi_\infty(\Xi_2)$. We conclude from \eref{uniformsigmabd} and \eref{pf3eqn1} that for $m=1,2,\cdots$,
\be\label{pf3eqn3}
\|\mathcal{E}(P)\|_{\XX_1,\mu_1^*,p}\le cn^q\|P\|_{\XX_2,\mu_2^*,p}, \qquad \|\mathcal{E}(P)-\sigma_{2^m}(P)\|_{\XX_1,\mu_1^*,p}\le cn^q2^{-m\beta}\|P\|_{\XX_2,\mu_2^*,p}, \qquad P\in\Pi_n(\Xi_2).
\ee

Next, let $f\in L^p(\XX_2,\mu_2^*)$. We find a sequence $\{P_m\}$ such that each $P_m\in\Pi_{2^m}(\Xi_2)$, and
$$
\|f-P_m\|_{\XX_2,\mu_2^*,p} \le cE_{2^m,p}(\Xi_2;f), \qquad m=0,1,\cdots.
$$
Then $P_{m+1}-P_m\in \Pi_{2^{m+1}}(\Xi_2)$, and a simple application of triangle inequality shows that
$$
\|P_{m+1}-P_m\|_{\XX_2,\mu_2^*,p}\le cE_{2^m,p}(\Xi_2;f).
$$
Therefore, if \eref{pseudobesovcond} holds, then
\eref{pf3eqn3} implies that
\be\label{pf3eqn4}
\|\mathcal{E}(P_0)\|_{\XX_1,\mu_1^*,p} +\sum_{m=0}^\infty \|\mathcal{E}(P_{m+1})-\mathcal{E}(P_m)\|_{\XX_1,\mu_1^*,p}
\le c\left\{\|P_0\|_{\XX_2,\mu_2^*,p} +\sum_{m=0}^\infty 2^{mq}E_{2^m,p}(\Xi_2;f)\right\}<\infty.
\ee
Consequently, there exists $\mathcal{E}(f)\in L^p(\XX_1,\mu_1^*)$ such that
$$
\mathcal{E}(f)=\mathcal{E}(P_0)+\sum_{m=0}^\infty \left(\mathcal{E}(P_{m+1})-\mathcal{E}(P_m)\right),
$$
with the series converging in the sense of $L^p(\XX_1,\mu_1^*)$. Moreover,
\be\label{pf3eqn5}
\|\mathcal{E}(f)-\mathcal{E}(P_n)\|_{\XX_1,\mu_1^*,p}  \le c\sum_{m=n}^\infty 2^{mq}E_{2^m,p}(\Xi_2;f).
\ee

Using this estimate, \eref{uniformsigmabd}, and \eref{pf3eqn3}, we conclude that
\begin{eqnarray*}
\|\mathcal{E}(f)-\sigma_{2^n}(f)\|_{\XX_1,\mu_1^*,p}&\le& \|\mathcal{E}(f)-\mathcal{E}(P_n)\|_{\XX_1,\mu_1^*,p} +\|\mathcal{E}(P_n)-\sigma_{2^n}(P_n)\|_{\XX_1,\mu_1^*,p} +\|\sigma_{2^n}(f-P_n)\|_{\XX_1,\mu_1^*,p}\\
&\le& c\sum_{m=n}^\infty 2^{mq}E_{2^m,p}(\Xi_2;f)+c2^{n(q-\beta)}\|P_n\|_{\XX_2,\mu_2^*,p}+c2^{nq}\|f-P_n\|_{\XX_2,\mu_2^*,p}\\
&\le& c\sum_{m=n}^\infty 2^{mq}E_{2^m,p}(\Xi_2;f)+c2^{n(q-\beta)}\|f\|_{\XX_2,\mu_2^*,p}.
\end{eqnarray*}
This proves \eref{twospacegoodapprox}, and in particular, the fact that $\mathcal{E}(f)=\lim_{n\to\infty}\sigma_{2^n}(f)$ in the sense of $L^p(\XX_1,\mu_1^*)$. 

We observe that the sum expression in \eref{jointsigmaopdef} involves non--zero terms only for $j, k$ for which $\ell_{j,k} <n$. If $\a\ell_{j,k}\ge \lambda_{1,j}$ for $j, k=0,1,\cdots$, then the only values of $j$ for which the summand in \eref{jointsigmaopdef} is non--zero  are those for which $\lambda_{1,j}<\a n$. Therefore, $\sigma_n(f)\in \Pi_{\a n}(\Xi_1)$.
Therefore, \eref{twospacegoodapprox} implies that
\be\label{pf3eqn6}
E_{\a 2^n,p}(\Xi_1;f)\le c\left\{\sum_{m=n}^\infty 2^{mq}E_{2^m,p}(\Xi_2;f) + 2^{n(q-\beta)}\|f\|_{\XX_2,\mu_2^*,p}\right\}.
\ee
If $f\in B_{p,\rho,\gamma}(\Xi_2)$, then by definition the sequence $\{2^{mq}E_{2^m,p}(\Xi_2,f)\}$ is in $\mathsf{b}_{\rho,\gamma-q}$. The discrete Hardy inequality \cite[Lemma~3.4, p.~27]{devlorbk} implies that the sequence
$$
\left\{\sum_{m=n}^\infty 2^{mq}E_{2^m,p}(\Xi_2;f)\right\}_{n=0}^\infty \in \mathsf{b}_{\rho,\gamma-q}.
$$
Since $\beta>\gamma>q$, then $\{2^{n(q-\beta)}\}_{n=0}^\infty \in \mathsf{b}_{\rho,\gamma-q}$ as well. Hence, \eref{pf3eqn6} leads to the fact that $\mathcal{E}(f)\in B_{p,\rho,\gamma-q}$. \qed

Next, we turn to the proof of Theorem~\ref{gentensortheo}. The analogue of Theorem~\ref{mainappltheo} in this case is the following theorem.

\begin{theorem}\label{tensopbdtheo}
Let $\Xi_1=(\XX_1,d_1,\mu_1^*,\{\lambda_{1,k}\}, \{\phi_{1,k}\})$ and
 $\Xi_2=(\XX_2,d_2,\mu_2^*,\{\lambda_{2,k}\}, \{\phi_{2,k}\})$ be two admissible systems as in Definiton~\ref{systemdef}, each of which satisfies the Gaussian upper bound condition with exponents $q_1$, $q_2$ respectively, where $\{\phi_{1,j}\}$ and $\{\phi_{2,k}\}$ are orthonormalized with respect to $\mu_1^*$, respectively, $\mu_2^*$. Let $\YY$ be a measurable (with respect to both $\mu_1^*$ and $\mu_2^*$) subset of $\XX_1\cap \XX_2$, $\nu^*$ be a probability measure on $\YY$, $\Gamma_{j,k}$ be defined as in  \eref{genconnectiondef}. We assume further that  regularity condition \eref{ballmeasurecond2}
 holds.
Let $1\le p\le\infty$, $f\in L^p(\XX_2,\mu_2^*)$, and $h$ be a sufficiently smooth low pass filter. Then
 \be\label{tensopbd}
\|\sigma_{n,\otimes}(\Xi_1,\Xi_2;h,f)\|_{\XX_1,\mu_1^*,p}\le cn^{q_1+(q_2-q_1)/p} \|f\|_{\XX_2,\mu_2^*, p}.
\ee
\end{theorem}

\begin{Proof}\ 
In this proof, we write for $j=1,2$, $x, y\in \XX_j$, $f\in L^1(\XX_j,\mu_j^*)$,
\be\label{pf10eqn1}
\Phi_{n,j}(x,y)=\sum_{k=0}^\infty h\left(\frac{\lambda_{j,k}}{n}\right)\phi_{j,k}(x)\overline{\phi_{j,k}(y)}, \qquad \sigma_{n,j}(f)(x)=\int_\XX \Phi_{n,j}(x,y)f(y)d\mu_j^*(y).
\ee
We observe from the definitions that
\be\label{pf10eqn2}
\Phi_{n, \otimes}(\Xi_1,\Xi_2;h, x_1, x_2)=\int_\YY \Phi_{n,1}(x_1,z)\overline{\Phi_{n,2}(x_2,z)}d\nu^*(z).
\ee
As in the proof of Lemma~\ref{fundalemma}, we see using \eref{ballmeasurecond2} that
$$
\int_{\XX_2}|\Phi_{n, \otimes}(\Xi_1,\Xi_2;h, x_1, x_2)|d\mu_2^*(x_2) \le \int_\YY |\Phi_{n,1}(x_1,z)|\int_{\XX_2}|\Phi_{n,2}(x_2,z)|d\mu_2^*(x_2)d\nu^*(z)\le cn^{q_1},
$$
and an analgous inequality holds with integrals on $\XX_1$ instead. Hence, using the Riesz--Thorin theorem as in the proof of Theorem~\ref{mainappltheo}, we arrive at \eref{tensopbd}.
\end{Proof}

\noindent
\textsc{Proof of Theorem~\ref{gentensortheo}.}
This theorem proved exactly in the same way as in the proof of Theorem~\ref{mainappltheo}, using Theorem~\ref{tensopbdtheo} in place of Theorem~\ref{mainappltheo}. The condition on the joint eigen--values in Theorem~\ref{extensiontheo} is superfluous in this context; it is clear that if $P\in\Pi_n(\Xi_2)$, then $\sigma_{n,\otimes}(\Xi_1,\Xi_2;h,P)\in\Pi_n(\Xi_1)$. We omit the remaining details of this proof. \qed

\subsection{Tauberian theorem}\label{taubsect}
We recall that if $\mu$ is an extended complex valued Borel measure on $\RR$, then its total variation measure is defined for a Borel set $B$ by
$$
|\mu|(B)=\sup\sum |\mu(B_k)|,
$$
where the sum is over a partition $\{B_k\}$ of $B$ comprising Borel sets, and the supremum is over all such partitions.

A measure $\mu$ on $\RR$ is called an even measure if $\mu((-u,u))=2\mu([0,u))$ for all $u>0$, and $\mu(\{0\})=0$. If $\mu$ is an extended complex valued measure on $[0,\infty)$, and $\mu(\{0\})=0$, we define a measure $\mu_e$ on $\RR$ by 
$$
\mu_e(B)=\mu\left(\{|x| : x\in B\}\right),
$$
and observe that $\mu_e$ is an even measure such that $\mu_e(B)=\mu(B)$ for $B\subset [0,\infty)$. In the sequel, we will assume that all measures on $[0,\infty)$ which do not associate a nonzero mass with the point $0$ are extended in this way, and will abuse the notation $\mu$ also to denote the measure $\mu_e$. In the sequel, the phrase ``measure on $\RR$'' will refer to an extended complex valued Borel measure having bounded total variation on compact intervals in $\RR$, and similarly for measures on $[0,\infty)$.

Our main Tauberian theorem is the following.

\begin{theorem}\label{maintaubertheo}
Let $\mu$ be an extended complex valued measure on $[0,\infty)$, and $\mu(\{0\})=0$. We assume that there exist $Q, r>0$, such that each of the following conditions are satisfied.
\begin{enumerate}
\item 
\be\label{muchristbd}
\tn\mu\tn_Q:=\sup_{u\in [0,\infty)}\frac{|\mu|([0,u))}{(u+2)^Q} <\infty,
\ee
\item There are constants $c, C >0$,  such that
\be\label{muheatgaussbd}
\left|\int_\RR \exp(-u^2t)d\mu(u)\right|\le c_1t^{-C}\exp(-r^2/t)\tn\mu\tn_Q, \qquad 0<t\le 1.
\ee 
\end{enumerate}
Let $H:[0,\infty)\to\RR$, $S>Q+1$ be an integer, and suppose that there exists a measure $H^{[S]}$ such that
\be\label{Hbvcondnew}
H(u)=\int_0^\infty (v^2-u^2)_+^{S}dH^{[S]}(v), \qquad u\in\RR,
\ee
and
\be\label{Hbvintbdnew}
V_{Q,S}(H)=\max\left(\int_0^\infty (v+2)^Qv^{2S}d|H^{[S]}|(v), \int_0^\infty (v+2)^Qv^Sd|H^{[S]}|(v)\right)<\infty.
\ee
Then for $n\ge 1$,
\be\label{genlockernest}
\left|\int_0^\infty H(u/n)d\mu(u)\right| \le c\frac{n^Q}{\max(1, (nr)^S)}V_{Q,S}(H)\tn\mu\tn_Q.
\ee
\end{theorem}
 
The proof of Theorem~\ref{maintaubertheo} will be given in a number of steps. First, we will show that the conclusions of the theorem are valid if they are valid when the function $H$ is replaced by a special function. The next step is to show that the theorem is valid if the distributional Fourier transform of $\mu$ is supported on $(-\infty,-2r]\cup [2r,\infty)$. The final step is to prove that the assumptions in Theorem~\ref{maintaubertheo} imply this property for the distributional Fourier transform.\\

\noindent
\textsc{Step 1.}\\

Let
\be\label{truncpowdef}
x_+^r=\left\{
\begin{array}{ll}
x^r, &\mbox{ if $x>0$,}\\
0, &\mbox{ if $x\le 0$.}
\end{array}\right.
\ee
The Bochner--Riesz summability method applied to $\mu$ is defined by
\be\label{BRmeasuredef}
\mathfrak{R}_{S;n}(\mu):=\int_\RR \left(1-\frac{u^2}{n^2}\right)_+^Sd\mu(u).
\ee

\begin{lemma}\label{BRequivlemma}
Let $S>Q>0$,  and $\mu$ be an even measure satisfying \eref{muchristbd} (and $\mu(\{0\})=0$).  If
\be\label{BRloccond}
\mathfrak{R}_{S;n}(\mu)\le  c\frac{(n+2)^Q}{\max(1, (nr)^S)}\tn\mu\tn_Q, \qquad n>0,
\ee
then \eref{genlockernest} holds for all $H$ satisfying \eref{Hbvcondnew}, \eref{Hbvintbdnew}.
\end{lemma}

\begin{Proof}\ 
Without loss of generality, we may assume that $\tn\mu\tn_Q=1$. Let $n\ge 1$. Assuming that a change of the order of integration below is valid,  we obtain from \eref{Hbvcondnew} that
\be\label{pf1eqn1}
\int_0^\infty H(u/n)d\mu(u)=\int_0^\infty\!\!\!\int_0^\infty (v^2-(u/n)^2)_+^SdH^{[S]}(v)d\mu(u)=(1/2)\int_0^\infty \mathfrak{R}_{S;nv}(\mu)v^{2S}dH^{[S]}(v).
\ee
We observe that for $n\ge 1$, $v\ge 0$,
$$
(nv+2)\le n(v+2).
$$
Therefore, the assumption \eref{BRloccond} implies that if $nr>1$ then
\bea\label{pf1eqn2}
\left|\int_0^\infty H(u/n)d\mu(u)\right|&\le& cn^Q\int_0^\infty \frac{(v+2)^Qv^{2S}d|H^{[S]}|(v)}{\max(1, (nrv)^S)}\nonumber\\
&=& cn^Q\left\{\int_0^{1/(nr)} (v+2)^Qv^{2S}d|H^{[S]}|(v) + (nr)^{-S}\int_{1/(nr)}^\infty (v+2)^Qv^Sd|H^{[S]}|(v)\right\}\nonumber\\
&=& \frac{cn^Q}{(nr)^S}\int_0^\infty (v+2)^Qv^Sd|H^{[S]}|(v).
\eea
If $nr\le 1$, then \eref{pf1eqn1} and \eref{BRloccond} together imply that
\be\label{pf1eqn3}
\left|\int_0^\infty H(u/n)d\mu(u)\right|\le cn^Q\int_0^\infty (v+2)^Qv^{2S}d|H^{[S]}|(v).
\ee
Using the condition \eref{Hbvintbdnew}, we deduce \eref{genlockernest} immediately from \eref{pf1eqn2}, \eref{pf1eqn3}.

It remains to verify that the interchange of order of integration in \eref{pf1eqn1} is justified. Using \eref{muchristbd} and \eref{Hbvintbdnew}, we observe that
\begin{eqnarray*}
\lefteqn{\int_0^\infty\!\!\!\int_0^\infty (v^2-(u/n)^2)_+^{S}d|\mu|(u)d|H^{[S]}|(v)= \int_0^\infty\!\!\!\int_0^{nv} (v^2-(u/n)^2)_+^{S}d|\mu|(u)d|H^{[S]}|(v)}\\
&\le& \int_0^\infty v^{2S}\int_0^{nv} d|\mu|(u) d|H^{[S]}|(v) \le n^Q\int_0^\infty (v+2)^Qv^{2S}d|H^{[S]}|(v)<\infty.
\end{eqnarray*}
This justifies the change of the order of integration, and completes the proof.
\end{Proof}

\begin{rem}\label{hintexistrem}
{\rm The verification for the interchange of order of integration shows in particular that if $V_{Q,S}(H)<\infty$, then $\int_0^\infty |H(u/n)|d|\mu|(u)<\infty$ for all $n\ge 1$. \qed
}
\end{rem}

\noindent
\textsc{Step 2.}\\

Here we work with the support of the distributional Fourier transform of $\mu$. We need first to recall a few facts about distributions and their Fourier transforms, as well as certain facts from measure theory, mainly to establish notation, and explain our abuse thereof.

A function $\phi :\RR\to\RR$ is called a rapidly decreasing function (or test function) if it is infinitely differentiable and $|x|^k\derf{\phi}{j}(x)\to 0$ as $|x|\to\infty$ for all non--negative integers $k, j$. The space of all test functions is denoted by $\mathcal{S}$, and is a locally convex linear space. Its (continuous) dual is denoted by $\mathcal{S}'$, and a member of this dual is called a (tempered) distribution. If $\phi\in \mathcal{S}$ its Fourier transform $\hat{\phi}$ and inverse Fourier transform $\tilde{\phi}$ are defined by
\be\label{testfourtransdef}
\hat{\phi}(x)=\frac{1}{2\pi}\int_\RR \phi(y)\exp(-ixy)dy, \qquad \tilde{\phi}(x) =\int_\RR \phi(y)\exp(ixy)dy,\qquad x\in\RR.
\ee
Both $\hat{\phi}$ and $\tilde{\phi}$ are test functions.
Likewise, if $\mu\in \mathcal{S}'$, its Fourier transform $\hat{\mu}$ and inverse Fourier transform $\tilde{\mu}$ are defined by
\be\label{distfourtransdef}
\hat{\mu}(\phi)=\mu(\hat{\phi}), \qquad \tilde{\mu}(\phi)=\mu(\tilde{\phi}), \qquad \phi\in\mathcal{S}.
\ee
Clearly, both $\hat{\mu}$ and $\tilde{\mu}$ are distributions.
One has the inversion formulas
$$
\widehat{\tilde{\phi}}=\widetilde{\hat{\phi}} =\phi, \qquad 
\widehat{\tilde{\mu}}=\widetilde{\hat{\mu}} =\mu, \qquad \phi\in \mathcal{S}, \ \mu\in \mathcal{S}'.
$$
We refer the reader to \cite{rudinreal, rudinfunctional} for further information regarding Fourier transforms, and their properties.

If $\mu$ is a signed, extended complex valued, even, Borel measure on $\RR$, we will identify it with the corresponding function on $\RR$, defined by $u\mapsto \mu([0, |u|))$, $u\in\RR$, and abuse the notation to denote both the measure and the function by $\mu$. The total variation measure  $|\mu|$ 
of $\mu$ is  also identified with the total variation function of the corresponding function on $\RR$, and with the same abuse of notation. If there exists $Q\ge 0$ such that $|\mu|((-u,u))\le cu^Q$, $u>0$, then $\mu$ also defines an element of $\mathcal{S}'$ by the formula
\be\label{measureasdist}
\mu(\phi)=\int_\RR\phi d\mu, \qquad |\mu|(\phi)=\int_\RR\phi d|\mu|, \qquad \phi\in \mathcal{S}.
\ee
One advantage of the identification between the measure, the function, and the distribution is that one can use such theorems from measure theory as Fubini's theorem, dominated convergence theorem, etc. in many cases involving integration of test functions with respect to extended complex valued measures satisfying the condition just described. In particular, it is often useful to integrate by parts in \eref{measureasdist} to observe that
\be\label{measureasdistbis}
\mu(\phi)=-\int_\RR\phi'(u)\mu(u)du, \qquad 
|\mu|(\phi)=-\int_\RR\phi'(u)|\mu|(u)du, \qquad \phi\in\mathcal{S}.
\ee
It is unfortunately customary to denote the measure, the function, and the distribution all by the same symbol, with the corresponding abuse of notations regarding the Fourier transforms, e.g.,
$$
\hat{\mu}(\phi)=\frac{1}{2\pi}\int_\RR \hat{\phi}d\mu, \qquad \phi\in\mathcal{S},
$$ 
and we will follow this convention, with clarifying remarks in case we feel a cause for confusion.  The support of a function $f :\RR\to\CC$, denoted by $\supp(f)$, is the closure of the set $\{x\in\RR : f(x)\not=0\}$, and $f$ is said to be supported on any super--set of this set. The support of a measure $\mu$ on $\RR$ that can be viewed as a distribution is the closure of the set $A\subset\RR$ such that $\int_\RR \phi d\mu =0$ for every $\phi\in \mathcal{S}$, supported on $\RR\setminus A$. 

We are now ready to state our auxiliary Tauberian theorem.

\begin{theorem}\label{auxtaubertheo}
Let $\mu$ be an extended complex valued, even measure supported on $\RR$ (with $\mu(\{0\})=0$). We assume that there exist $Q, r>0$, such that \eref{muchristbd} holds, and 
the (distributional) Fourier transform $\hat{\mu}$ of $\mu$ is supported on $(-\infty, -2r]\cup[2r,\infty)$.
Then the conclusion of Therorem~\ref{maintaubertheo} holds.
\end{theorem}

In order to prove this theorem, we will show that the assumptions imply \eref{BRloccond}. First, we need a technical lemma.

\begin{lemma}\label{bandlimitapproxlemma}
Let $S\ge 2$ be an integer. For every $Y>1$, there exists an even function $H_Y\in\mathcal{S}$ such that
\be\label{bandlimit}
\widehat{H_Y}(x)=0, \qquad |x|\ge Y,
\ee
\be\label{brapprox}
\left|\frac{d^k}{du^k}(1-u^2)_+^{S}-H_Y^{(k)}(u)\right|\le cY^{k-S}, \qquad u\in\RR,\ k=0,1,\cdots, S-1,
\ee
and for any integer $N\ge 1$,
\be\label{brlocapprox}
|H_Y^{(k)}(u)|\le c(N)(Y|u|)^{-N}, \qquad |u|\ge 2,\ k=0,1,\cdots, S-1.
\ee
\end{lemma}

\begin{Proof}\ 
The proof is by a direct construction. In this proof, we denote $F(u)=(1-u^2)_+^{S}$. Since $S\ge 2$, $F$ has a first derivative having a bounded total variation on $\RR$. Therefore,  the Fourier inversion formula holds for $F$. The Fourier transform of $F$ is given in terms of the Bessel functions \cite[p.~53]{watson} by
\bea\label{besselfn}
\hat{F}(x)&=&\frac{S!2^{S-1/2}}{\sqrt{\pi}}\frac{J_{S+1/2}(x)}{x^{S+1/2}}\nonumber\\
&=&\frac{S!2^{S-1}}{\pi}\frac{1}{x^{S+1}}\left\{e^{ix}\sum_{k=0}^{S}\frac{i^{k-S-1}(S+k)!}{k!(S-k)! (2x)^k}+e^{-ix}\sum_{k=0}^{S}\frac{(-i)^{k-S-1}(S+k)!}{k!(S-k)! (2x)^k}\right\}.
\eea
Since $F$ is supported on $[-1,1]$, $\hat{F}$ is an entire function of finite exponential type $1$. In this proof, let $h$ be an infinitely differentiable low pass filter. (Necessarily, $h$ is a test function.)  We define $H_Y$ by
\be\label{pf2eqn1}
\widehat{H_Y}(x)=\hat{F}(x)h(x/Y), \qquad x\in\RR,
\ee 
or equivalently,
\be\label{pf2eqn2}
H_Y(u)=\frac{Y}{2\pi}\int_\RR F(v)\tilde{h}(Y(u-v))dv=\frac{Y}{2\pi}\int_{-1}^1 F(v)\tilde{h}(Y(u-v))dv, \qquad u\in\RR.
\ee
It is not difficult to verify from the definitions that $\widehat{H_Y}\in\mathcal{S}$, and hence, $H_Y\in\mathcal{S}$ also . The equation \eref{bandlimit} is also clear from the definition of $\widehat{H_Y}$. 

In this proof, let $g(u)=h(u)-h(2u)$, $u\in\RR$, and for every integer $n\ge 1$,  $G_{n,Y}$ be defined by
\be\label{pf2eqn3}
\widehat{G_{n,Y}}(x)=\hat{F}(x)g(x/(2^nY)),\qquad G_{n,Y}(u)=\frac{2^nY}{2\pi}\int_\RR F(v)\tilde{g}(2^nY(u-v))dv, \qquad x, u\in\RR.
\ee
For every integer $n\ge 1$, $g(x/(2^nY))=0$ if $|x|\le 2^{n-2}Y
$ or $|x|\ge 2^nY$. Therefore, the second equation in \eref{besselfn} shows that for $0\le k\le S-1$,
\be\label{pf2eqn5}
\left| \widehat{G_{n,Y}^{(k)}}(x)\right|=|x|^{k}|\widehat{G_{n,Y}}(x)|\le c(2^nY)^{k-S-1}, \qquad x\in\RR.
\ee
Since $\widehat{G_{n,Y}^{(k)}}$ is supported on $[-2^nY, -2^{n-2}Y]\cup [2^{n-2}Y, 2^nY]$, the Fourier inversion formula implies that
$$
\left|G_{n,Y}^{(k)}(u)\right|\le \left|\int_{2^{n-2}Y\le |x|\le 2^nY} e^{iux}\widehat{G_{n,Y}^{(k)}}(x)dx\right| \le c(2^nY)^{k-S}, \qquad u\in\RR, \ k=0,1,\cdots, S-1.
$$
Consequently,
\be\label{pf2eqn6}
\sum_{n=1}^\infty \left|G_{n,Y}^{(k)}(u)\right|\le cY^{k-S}, \qquad u\in\RR, \ k=0,1,\cdots, S-1,
\ee
and hence,
\be\label{pf2eqn4}
F^{(k)}(u)=H_Y^{(k)}(u)+\sum_{n=1}^\infty G_{n,Y}^{(k)}(u), \qquad  u\in\RR, \ k=0,1,\cdots, S-1,
\ee
where the series converges uniformly and absolutely on $\RR$. Moreover, \eref{brapprox} is clear from \eref{pf2eqn4} and \eref{pf2eqn6}.

Next, we observe that for $|u|\ge 2$, $|v|\le 1$, $|u-v|\ge (1/2)|u|$. 
Since $\tilde{h}$ is a test function, \eref{pf2eqn2} implies that for $|u|\ge 2$,
$$
|H_Y^{(k)}(u)|\le cY^{k+1}\int_{-1}^1 |F(v)||\tilde{h}^{(k)}
(Y(u-v))|dv \le c(N)Y^{k+1}(Y|u|)^{-k-1-N}\int_{-1}^1 |F(v)|dv\le c(N)(Y|u|)^{-N}.
$$
\end{Proof}

We are now ready to prove Theorem~\ref{auxtaubertheo}.

\noindent
\textsc{Proof of Theorem~\ref{auxtaubertheo}.}  In this proof, we may assume without loss of generality that $\tn\mu\tn_Q=1$. We will prove \eref{BRloccond}. This estimate is clear immediately from \eref{muchristbd} if $nr\le 1$. In this proof, we assume $nr>1$, let $Y=2nr$, and $H_Y$ denote the
function as in Lemma~\ref{bandlimitapproxlemma}. Then $\widehat{H_Y(\cdot/n)}$ is supported on $(-2r,2r)$. Since $\hat{\mu}$ is supported on $(-\infty, -2r]\cup [2r,\infty)$,
\be\label{pf8eqn1}
\int_\RR H_Y(u/n)d\mu(u)=0.
\ee
Using \eref{muchristbd}, and \eref{brapprox} with $k=0$, we deduce that
\be\label{pf8eqn2}
\left|\int_{-2n}^{2n} \left\{\left(1-\frac{u^2}{n^2}\right)^S_+ - H_Y(u/n)\right\}d\mu(u)\right| \le c\frac{(n+2)^Q}{(nr)^S}.
\ee
If $|u|\ge 2n$, then
$$
\left(1-\frac{u^2}{n^2}\right)^S_+ - H_Y(u/n)=-H_Y(u/n).
$$
Using integration by parts and \eref{brlocapprox} with $N=S>Q+1$, we obtain
\bea\label{pf8eqn3}
\lefteqn{\left|\int_{|u|\ge 2n} \left\{\left(1-\frac{u^2}{n^2}\right)^S_+ - H_Y(u/n)\right\}d\mu(u)\right|}\nonumber\\
&=&\left|\int_{|u|\ge 2n} H_Y(u/n)d\mu(u)\right|\le |\mu|(2n)H_Y(2)+\frac{1}{n}\left|\int_{|u|\ge 2n} H_Y'(u/n)\mu(u)du\right|\nonumber\\
&\le &  c(n+2)^Q(nr)^{-S} +\frac{1}{nr^S}\int_{|u|\ge 2n}u^{Q-S}du\nonumber\\
&\le& cn^Q(nr)^{-S}.
\eea
The estimate \eref{BRloccond} follows from \eref{pf8eqn1}, \eref{pf8eqn2}, \eref{pf8eqn3}. In turn, Lemma~\ref{BRequivlemma} now implies Theorem~\ref{auxtaubertheo}. \qed

\noindent
\textsc{Step 3.}\\

We will prove that the assumptions on the measure in Theorem~\ref{maintaubertheo} imply those in Theorem~\ref{auxtaubertheo}; i.e., we prove the following theorem using some of the ideas in \cite{sikora2004riesz, frankbern}.

\begin{theorem}\label{sikoratheo}
Let $\mu$ be an even measure satisfying \eref{muchristbd}. If there exist   $r, C>0$ such that \eref{muheatgaussbd} is satisfied.
Then the  (distributional) Fourier transform $\hat{\mu}$ of $\mu$ is supported on $(-\infty, -2r]\cup[2r,\infty)$.
\end{theorem}

The proof of this theorem is fairly complicated, and is organized in the form of three preparatory lemmas.
\begin{lemma}\label{bvapproxlemma}
Let $\mu$ be an even measure satisfying \eref{muchristbd}, and $R\ge 4$. If $\phi$ is a test function, and $h$ is an infinitely differentiable low pass filter, then for every integer $N\ge 2$, 
\be\label{bvapprox}
\left|\int_\RR \phi(u)h(u/R)d\mu(u)-\int_\RR \phi(u)d\mu\right| \le \frac{2^{Q+N+1}}{NR^N}\left\{\max_{|x|\ge R/2}|\phi(u)u^{Q+N}| +\max_{|x|\ge R/2}|\phi'(u)u^{Q+N+1}|\right\}\tn\mu\tn_Q\max_{x\in\RR}|h'(x)|.
\ee
\end{lemma}

\begin{Proof}\ 
Without loss of generality, we may assume that $\tn\mu\tn_Q=1$. Since $h$ is a low pass filter, $1-h(u/R)=0$ and $h'(u/R)=0$ if $|u|< R/2$. If $|u|\ge R/2\ge 2$, \eref{muchristbd} implies that $|\mu|(u)\le 2^Q|u|^Q$. Therefore, using integration by parts and \eref{muchristbd}, we obtain
\bea\label{pf5eqn1}
\left|\int_\RR \phi(u)h(u/R)d\mu(u)-\int_\RR \phi(u)d\mu(u)\right|&=&\left|\int_\RR \left\{\phi'(u)(1-h(u/R))-\frac{1}{R}\phi(u)h'(u/R)\right\}\mu(u)du\right|\nonumber\\
&\le& 2^Q\int_{|u|\ge R/2}|\phi'(u)u^{Q+N+1}|\frac{du}{|u|^{N+1}} \nonumber\\
&&\qquad+\max_{x\in \RR}|h'(x)|\frac{2^Q}{R}\int_{|u|\ge R/2}|\phi(u)u^{Q+N}|\frac{du}{|u|^{N}}.
\eea
We observe that since $h$ is low pass filter, the mean value theorem implies that $\max_{x\in \RR}|h'(x)|\ge 2$. The estimate \eref{bvapprox} is now easy to deduce.
\end{Proof}

\begin{lemma}\label{kanallemma}
 Let $\mu$ be an even measure satisfying \eref{muchristbd}, and   $\phi : \RR \to\RR$ be an infinitely differentiable function supported on $[0,b^2]$ for some $b>0$. Let $\psi :\RR\to\RR$ be defined by $\psi(u)=\phi(u^2)$.\\
 {\rm (a)} The function $\psi\in\mathcal{S}$.\\ 
{\rm (b)}  The function $\hat{\phi}$ can be extended to $\CC$ as an entire function of finite exponential type $\le b^2$, and for any integer $m\ge 0$,
\be\label{phibd}
|\hat{\phi}(z)|\le \frac{c\exp(b^2|\Im  z|)}{|z|^m}\int_\RR|\derf{\phi}{m}(u)|du, \qquad z\in\CC\setminus\{0\}.
\ee
{\rm (c)} If
\be\label{kanaldef}
F(z)=\int_0^\infty \exp\left(-\frac{u^2}{4iz}\right)d\mu(u), \qquad z=x+iy,\ y<0,
\ee
 then
\be\label{psimutokanal}
\hat{\mu}(\psi)=\int_\RR \hat{\phi}(-x-iy)(i(x+iy)/\pi)^{-1/2}F(x+iy)dx, \qquad y<0,
\ee
where the principal branch of the square root is understood.
\end{lemma}

\begin{Proof}\ 
Part (a) is easy to verify using the chain rule of differentiation several times. 

Defining
$$
\hat{\phi}(z)=\frac{1}{2\pi}\int_0^{b^2} \exp(-iuz)\phi(u)du, \qquad z\in\CC,
$$
it is clear that $\hat{\phi}$ extends the function $\hat{\phi}$ on $\RR$ to $\CC$. The fact that this extension is an entire function of finite exponential type $\le b^2$ and \eref{phibd} are elementary to verify. This proves part (b).

In view of \eref{phibd}, we may apply the Fourier inversion formula to deduce that
\be\label{pf6eqn1}
\phi(u)=\exp(uy)\int_\RR \hat{\phi}(-x-iy)\exp(-ixu)dx, \qquad u, y\in\RR.
\ee

Without loss of generality, we may assume that $\tn\mu\tn_Q=1$.
In this proof only, let $h$ be an infinitely differentiable low pass filter as in Lemma~\ref{bvapproxlemma} and for $u\in\RR$, $d\mu_R(u)=h(u/R)d\mu(u)$, $R\ge 4$. Then $\mu_R$ is a function with bounded total variation on $\RR$, $|\mu_R|(\RR)\le cR^Q$, $R\ge 2$. Hence, $\widehat{\mu_R}$ is an everywhere defined, even, and bounded function on $\RR$, with $|\widehat{\mu_R}(u)|\le cR^Q$, $u\in\RR$. In this proof, we write
$$
g_R(u)=u^{-1/2}\widehat{\mu_R}(\sqrt{u}), \qquad u\in [0,\infty), \ R\ge 4.
$$
Then
$$
\int_0^\infty \phi(u)g_R(u)du =2\int_0^\infty \phi(u^2)\widehat{\mu_R}(u)du=\int_\RR\psi(u)\widehat{\mu_R}(u)du =\int_\RR \hat{\psi}(u)d\mu_R(u),
$$
and Lemma~\ref{bvapproxlemma} shows that
\be\label{pf6eqn2}
\hat{\mu}(\psi)=\int_\RR \hat{\psi}(u)d\mu(u)=\lim_{R\to\infty} \int_0^\infty \phi(u)g_R(u)du.
\ee
Next, we recall \cite[Lemma~7.6]{rudinfunctional} that
$$
\exp(-x^2/2)=\frac{1}{\sqrt{2\pi}}\int_\RR \exp(-v^2/2)\exp(-ivx)dv.
$$
Let $t>0$. Letting $x=u/\sqrt{2t}$, and $v=w\sqrt{2t}$, this leads  to
$$
\exp(-u^2/(4t))=\sqrt{t/\pi}\int_\RR \exp(-tw^2)\exp(-iwu)dw, \qquad u\in\RR.
$$
Since $\mu_R$ is a measure of bounded total variation on $\RR$, we may use Fubini's theorem to conclude that
\be\label{pf6eqn3}
\int_\RR \exp(-u^2/(4t))d\mu_R(u)=\sqrt{t/\pi}\int_\RR \exp(-w^2t)\widehat{\mu_R}(w)dw=\sqrt{t/\pi}\int_0^\infty \exp(-tv)g_R(v)dv,
\ee
and hence, by analytic continuation that
\be\label{pf6eqn4}
F_R(z):=\int_\RR \exp(-u^2/(4iz))d\mu_R(u)=\sqrt{iz/\pi}\int_0^\infty \exp(-izv)g_R(v)dv, \qquad z\in\CC,\ \Im z <0.
\ee
In view of \eref{phibd}, we see that $\exp(uy)\hat{\phi}(-x-iy)g_R(u)$ is absolutely integrable with respect to $dxdu$ on $\RR\times [0,\infty)$ if $y<0$. Using \eref{pf6eqn1}, \eref{pf6eqn4}, we deduce that for $y<0$, $R\ge 2$,
\bea\label{pf6eqn5}
\int_0^\infty \phi(u)g_R(u)du &=& \int_0^\infty g_R(u)\left\{\exp(uy)\int_\RR \hat{\phi}(-x-iy)\exp(-ixu)dx\right\}du\nonumber\\
&=&\int_\RR \hat{\phi}(-x-iy)\left\{\int_0^\infty\exp(-i(x+iy)u)g_R(u)du\right\}dx\nonumber\\
&=&\int_\RR \hat{\phi}(-x-iy)(i(x+iy)/\pi)^{-1/2}F_R(x+iy)dx.
\eea
If $z=x+iy$, $y<0$, then it is elementary calculus to check that for every $L>0$, $|u^L\exp(-u^2/(4iz))|\le c(L)|z|^L|y|^{-L/2}$, and hence, \eref{phibd} implies that for such $z$, 
$$
|\hat{\phi}(-x-iy)(i(x+iy)/\pi)^{-1/2}F_R(x+iy)|\le c(L, y)|z|^{-2}.
$$
Therefore, using Lemma~\ref{bvapproxlemma}, the dominated convergence theorem, and \eref{pf6eqn2}, we arrive at \eref{psimutokanal} by letting $R\to\infty$ in \eref{pf6eqn5}.
This proves part (c).
\end{Proof}

\begin{lemma}\label{ktoanalcriticallemma}
Let $\mu$ be an even measure satisfying \eref{muchristbd} and \eref{muheatgaussbd}, and $F$ be defined by \eref{kanaldef}. Then
\be\label{ktoanalcriticalest}
|F(x+iy)|\le c|z|^{A}\exp(4r^2)\exp(4r^2y), \qquad y\le -1, \ x\in\RR,
\ee
where $A$ is the least integer $\ge \max(C,Q)$.
\end{lemma}

The proof of this lemma requires the Phragm\'en--Lindel\"of theorem \cite[Vol.~II, Theorem~7.5]{markushevich}, the required part of which we reproduce for the convenience of the reader.
\begin{theorem}\label{phragmentheo}
Let $D$ be the interior of an angle of $\a\pi$, $0\le \a\le 2$, with boundary $\Gamma$, and $f$ be analytic on $D$. With 
$$
M(f,v)=\sup_{|z|=v, z\in D}|f(z)|,
$$
let
\be\label{phragboundarycond}
\limsup_{z\to\xi}|f(z)|\le C<\infty, \qquad \xi\in\Gamma\cap \CC,
\ee
and
\be\label{phragdomaingrowthcond}
\liminf_{v\to\infty,\ v\in D}\frac{\log\log M(f,v)}{\log v} <1/\a,
\ee
then
\be\label{phargconclude}
|f(z)|\le C, \qquad z\in D.
\ee
\end{theorem}

\noindent
\textsc{Proof of Lemma~\ref{ktoanalcriticallemma}.} We will use the Phragm\'en--Lindel\"of  theorem with
$$
G(z)=\exp(4r^2iz)z^{-A}F(z), \qquad z\in\CC, \ \Im z<0,
$$
in place of $f$.
We note that $G$ is analytic in the domain as indicated. The estimate \eref{muheatgaussbd} implies that
$$
|F(-it)|=\left|\int_0^\infty \exp(-u^2/(4t))d\mu(u)\right|\le ct^C\exp(-4r^2i(-it)), \qquad t \ge 1/4.
$$
In particular, 
\be\label{pf7eqn1}
|G(z)|\le c, \qquad z=-it,\ t\ge 1/4.
\ee

Also, using \eref{measureasdistbis} and \eref{muchristbd}, we deduce that for $z=x-iy$, $x\in\RR$, $y>0$,
\bea\label{pf7eqn2}
|F(x-iy)|&\le& \int_0^\infty \left|\exp\left(-\frac{u^2}{4i(x-iy)}\right)\right|d|\mu|(u)=\int_0^\infty \exp(-u^2y/(4|z|^2))d|\mu|(u)\nonumber\\
&=&\frac{y}{2|z|^2}\int_0^\infty u\exp(-u^2y/(4|z|^2))|\mu|(u)du\nonumber\\
&\le& \frac{cy}{|z|^2}\int_0^\infty u^{Q+1}\exp(-u^2y/(4|z|^2))du\nonumber\\
&=& \frac{c|z|^Q}{y^{Q/2}}\int_0^\infty u^{Q/2}e^{-u}du.
\eea
Using this estimate with $z=x-i$, we obtain that
\be\label{pf7eqn3}
|G(z)|=\left|\exp(4r^2iz)z^{-A}F(z)\right|\le c\exp(4r^2), \qquad z=x-i,\ x\in\RR.
\ee
The estimate \eref{pf7eqn2} shows further that
 the condition \eref{phragdomaingrowthcond} is satisfied if $z\in\CC$, $\Im z \le -1$. Thus, we may apply the Phragm\'en--Lindel\"of  theorem with $G$, once taking  $D$ to be the sector bounded by $\{z=x-i,\ x\ge 0\}$ and $\{z=-it,\ t\ge 1\}$ and then with the  sector bounded by $\{z=x-i,\ x\le 0\}$ and $\{z=-it,\ t\ge 1\}$ to deduce \eref{ktoanalcriticalest} from \eref{pf7eqn1} and \eref{pf7eqn3}.\qed

\noindent
\textsc{Proof of Theorem~\ref{sikoratheo}.} Let $0<b<2r$, $\e>0$ and $\psi\in\mathcal{S}$ be an even function supported on $[-b,-\e]\cup[\e,b]$. Then there is an infinitely differentiable function $\phi$ supported on $[\e^2,b^2]$ such that $\psi(u)=\phi(u^2)$. With $F$ as in \eref{kanaldef}, Lemma~\ref{kanallemma} and Lemma~\ref{ktoanalcriticallemma} imply that for $y< -1$,
$$
|\hat{\mu}(\psi)| =\left|\int_\RR \hat{\phi}(-x-iy
)(i(x+iy)/\pi)^{-1/2}F(x+iy)dx\right|\le c\exp(4r^2)\exp(-(4r^2-b^2)|y|).
$$
Since $0<b<2r$, the right hand side of this inequality tends to $0$ as $y\to -\infty$. Hence $\hat{\mu}(\psi)=0$. Since $\e$, $b$, and $\psi$ were arbitrary, this implies that $\hat{\mu}$ is supported on $[-\infty,-2r]\cup [2r,\infty)$. \qed

Finally, we are in a position to put together the results in the three steps above to complete the proof of Theorem~\ref{maintaubertheo}.

\textsc{Proof of Theorem~\ref{maintaubertheo}}. In view of Theorem~\ref{sikoratheo}, the condition \eref{muheatgaussbd} implies that the support of $\hat{\mu}$ is contained in $(-\infty,-2r]\cup [2r,\infty)$. Theorem~\ref{maintaubertheo} follows immediately from Theorem~\ref{auxtaubertheo}. \qed

\bibliographystyle{plain}

\begin{thebibliography}{10}

\bibitem{askeyspecial}
R.~Askey.
\newblock {\em Orthogonal polynomials and special functions}.
\newblock SIAM, 1975.

\bibitem{niyogi1}
M.~Belkin, I.~Matveeva, and P.~Niyogi.
\newblock Regularization and semi-supervised learning on large graphs.
\newblock In {\em Learning theory}, pages 624--638. Springer, 2004.

\bibitem{belkin2003laplacian}
M.~Belkin and P.~Niyogi.
\newblock Laplacian eigenmaps for dimensionality reduction and data
  representation.
\newblock {\em Neural computation}, 15(6):1373--1396, 2003.

\bibitem{niyogi2}
M.~Belkin and P.~Niyogi.
\newblock Semi-supervised learning on {R}iemannian manifolds.
\newblock {\em Machine learning}, 56(1-3):209--239, 2004.

\bibitem{belkinfound}
M.~Belkin and P.~Niyogi.
\newblock Towards a theoretical foundation for {L}aplacian-based manifold
  methods.
\newblock {\em Journal of Computer and System Sciences}, 74(8):1289--1308,
  2008.

\bibitem{arjuna1}
A.~L. Bertozzi and A.~Flenner.
\newblock Diffuse interface models on graphs for classification of high
  dimensional data.
\newblock {\em Multiscale Modeling \& Simulation}, 10(3):1090--1118, 2012.

\bibitem{bochner_chadra_book}
S.~Bochner and K.~Chandrasekharan.
\newblock {\em Fourier transforms}.
\newblock Number~19 in Annals of mathematics studies. Princeton University
  Press, 1949.

\bibitem{achaspissue}
C.~K. Chui and D.~L. Donoho.
\newblock Special issue: Diffusion maps and wavelets.
\newblock {\em Appl. and Comput. Harm. Anal.}, 21(1), 2006.

\bibitem{chuiinterp}
C.~K. Chui and H.~N. Mhaskar.
\newblock Smooth function extension based on high dimensional unstructured
  data.
\newblock {\em Mathematics of Computation}, 83(290):2865--2891, 2014.

\bibitem{chuihyper}
C.~K. Chui and J.~Wang.
\newblock Dimensionality reduction of hyperspectral imagery data for feature
  classification.
\newblock In {\em Handbook of Geomathematics}, pages 1005--1047. Springer,
  2010.

\bibitem{chuiwang2010}
C.~K. Chui and J.~Wang.
\newblock Randomized anisotropic transform for nonlinear dimensionality
  reduction.
\newblock {\em GEM-International Journal on Geomathematics}, 1(1):23--50, 2010.

\bibitem{chuidimred2015}
C.~K. Chui and J.~Wang.
\newblock Nonlinear methods for dimensionality reduction.
\newblock In {\em Handbook of Geomathematics}, pages 1--46. Springer, 2015.

\bibitem{chung_directed_laplacian}
F.~Chung.
\newblock Laplacians and the cheeger inequality for directed graphs.
\newblock {\em Annals of Combinatorics}, 9(1):1--19, 2005.

\bibitem{coifman2013bi}
R.~R. Coifman and M.~J. Hirn.
\newblock Bi-stochastic kernels via asymmetric affinity functions.
\newblock {\em Applied and Computational Harmonic Analysis}, 35(1):177--180,
  2013.

\bibitem{coifmanhirn}
R.~R. Coifman and M.~J. Hirn.
\newblock Diffusion maps for changing data.
\newblock {\em Applied and Computational Harmonic Analysis}, 36(1):79--107,
  2014.

\bibitem{coifmanlafondiffusion}
R.~R. Coifman and S.~Lafon.
\newblock Diffusion maps.
\newblock {\em Applied and computational harmonic analysis}, 21(1):5--30, 2006.

\bibitem{coifmanmauro2006}
R.~R. Coifman and M.~Maggioni.
\newblock Diffusion wavelets.
\newblock {\em Applied and Computational Harmonic Analysis}, 21(1):53--94,
  2006.

\bibitem{daubbook}
I.~Daubechies.
\newblock {\em Ten lectures on wavelets}, volume~61.
\newblock SIAM, 1992.

\bibitem{david2003hessian}
L.~D. David and G.~Carrie.
\newblock Hessian eigenmaps: new locally linear embedding techniques for high
  dimensional data, tr2003-08, dept. of statistics, 2003.

\bibitem{davies1990heat}
E.~B. Davies.
\newblock {\em Heat kernels and spectral theory}, volume~92.
\newblock Cambridge University Press, 1990.

\bibitem{devlorbk}
R.~A. DeVore and G.~G. Lorentz.
\newblock {\em Constructive approximation}, volume 303.
\newblock Springer Science \& Business Media, 1993.

\bibitem{donoho2005image}
D.~L. Donoho and C.~Grimes.
\newblock Image manifolds which are isometric to euclidean space.
\newblock {\em Journal of mathematical imaging and vision}, 23(1):5--24, 2005.

\bibitem{donoho2002multiscale}
D.~L. Donoho, O.~Levi, J.-L. Starck, and V.~Martinez.
\newblock Multiscale geometric analysis for 3d catalogs.
\newblock In {\em Astronomical Telescopes and Instrumentation}, pages 101--111.
  International Society for Optics and Photonics, 2002.

\bibitem{compbio}
M.~Ehler, F.~Filbir, and H.~N. Mhaskar.
\newblock Locally learning biomedical data using diffusion frames.
\newblock {\em Journal of Computational Biology}, 19(11):1251--1264, 2012.

\bibitem{frankbern}
F.~Filbir and H.~N. Mhaskar.
\newblock A quadrature formula for diffusion polynomials corresponding to a
  generalized heat kernel.
\newblock {\em Journal of Fourier Analysis and Applications}, 16(5):629--657,
  2010.

\bibitem{modlpmz}
F.~Filbir and H.~N. Mhaskar.
\newblock Marcinkiewicz--{Z}ygmund measures on manifolds.
\newblock {\em Journal of Complexity}, 27(6):568--596, 2011.

\bibitem{filbirmadychacha}
F. Filbir, R. Hielscher, and W.R.~Madych.
\newblock Reconstruction from circular and spherical mean data.
\newblock {\em Applied and Computational Harmonic Analysis}, 29(1):111--120,
  2010.

\bibitem{friedman2004wave}
J.~Friedman and J.-P. Tillich.
\newblock Wave equations for graphs and the edge-based laplacian.
\newblock {\em Pacific Journal of Mathematics}, 216(2):229--266, 2004.

\bibitem{gidelew2014topics}
G.~A. Gidelew.
\newblock {\em Topics in harmonic analysis on combinatorial graphs}.
\newblock PhD thesis, Drexel University, 2014.

\bibitem{girosi1990networks}
F.~Girosi and T.~Poggio.
\newblock Networks and the best approximation property.
\newblock {\em Biological cybernetics}, 63(3):169--176, 1990.

\bibitem{gohbergkrein}
I.~Gohberg and M.~G. Krein.
\newblock {\em Introduction to the theory of linear nonselfadjoint operators},
  volume~18.
\newblock American Mathematical Soc., 1969.

\bibitem{grigoryan1995upper}
A.~Grigoryan.
\newblock Upper bounds of derivatives of the heat kernel on an arbitrary
  complete manifold.
\newblock {\em Journal of Functional Analysis}, 127(2):363--389, 1995.

\bibitem{grigor1997gaussian}
A.~Grigor'yan.
\newblock Gaussian upper bounds for the heat kernel on arbitrary manifolds.
\newblock {\em history}, 4(2exp):4t, 1997.

\bibitem{grigor1999estimates}
A.~Grigor’yan.
\newblock Estimates of heat kernels on {R}iemannian manifolds.
\newblock {\em London Math. Soc. Lecture Note Ser}, 273:140--225, 1999.

\bibitem{grigor2014heat}
A.~Grigor’yan, J.~Hu, and K.-S. Lau.
\newblock Heat kernels on metric measure spaces.
\newblock In {\em Geometry and Analysis of Fractals}, pages 147--207. Springer,
  2014.

\bibitem{halmosbdintop}
P.~R. Halmos and V.~S. Sunder.
\newblock {\em Bounded integral operators on $L^2$ spaces}, volume~96.
\newblock Springer Science \& Business Media, 2012.

\bibitem{hammond}
D.~K. Hammond, P.~Vandergheynst, and R.~Gribonval.
\newblock Wavelets on graphs via spectral graph theory.
\newblock {\em Applied and Computational Harmonic Analysis}, 30(2):129--150,
  2011.

\bibitem{niyogiface}
X.~He, S.~Yan, Y.~Hu, P.~Niyogi, and H.-J. Zhang.
\newblock Face recognition using {L}aplacianfaces.
\newblock {\em Pattern Analysis and Machine Intelligence, IEEE Transactions
  on}, 27(3):328--340, 2005.

\bibitem{horn_johnson_book}
R.~A. Horn and C.~R. Johnson.
\newblock {\em Matrix analysis}.
\newblock Cambridge university press, 2012.

\bibitem{jones2008parameter}
P.~W. Jones, M.~Maggioni, and R.~Schul.
\newblock Manifold parametrizations by eigenfunctions of the {L}aplacian and
  heat kernels.
\newblock {\em Proceedings of the National Academy of Sciences},
  105(6):1803--1808, 2008.

\bibitem{jones2010universal}
P.~W. Jones, M.~Maggioni, and R.~Schul.
\newblock Universal local parametrizations via heat kernels and eigenfunctions
  of the {L}aplacian.
\newblock {\em Ann. Acad. Sci. Fenn. Math.}, 35:131--174, 2010.

\bibitem{kim2014multi}
W.~H. Kim, V.~Singh, M.~K. Chung, C.~Hinrichs, D.~Pachauri, O.~C. Okonkwo,
  S.~C. Johnson, Alzheimer's Disease~Neuroimaging Initiative, et~al.
\newblock Multi-resolutional shape features via non-euclidean wavelets:
  Applications to statistical analysis of cortical thickness.
\newblock {\em NeuroImage}, 93:107--123, 2014.

\bibitem{kordyukov1991p}
Y.~A Kordyukov.
\newblock ${L}^p$--theory of elliptic differential operators on manifolds of
  bounded geometry.
\newblock {\em Acta Applicandae Mathematica}, 23(3):223--260, 1991.

\bibitem{korevaar_taub_book}
J.~Korevaar.
\newblock {\em Tauberian theory: a century of developments}, volume 329.
\newblock Springer Science \& Business Media, 2004.

\bibitem{lafon}
S.~S. Lafon.
\newblock {\em Diffusion maps and geometric harmonics}.
\newblock PhD thesis, Yale University, 2004.

\bibitem{ageface2011}
Z.~Li, U.~Park, and A.~K. Jain.
\newblock A discriminative model for age invariant face recognition.
\newblock {\em Information Forensics and Security, IEEE Transactions on},
  6(3):1028--1037, 2011.

\bibitem{marockner1992}
Z.-M. Ma and M.~R{\"o}ckner.
\newblock {\em Introduction to the theory of (non-symmetric) Dirichlet forms}.
\newblock Springer, 1992.

\bibitem{mauropap}
M.~Maggioni and H.~N. Mhaskar.
\newblock Diffusion polynomial frames on metric measure spaces.
\newblock {\em Applied and Computational Harmonic Analysis}, 24(3):329--353,
  2008.

\bibitem{markushevich}
A.~I. Markushevich.
\newblock {\em Theory of functions of a complex variable}, volume 296.
\newblock American Mathematical Soc., 2005.

\bibitem{eignet}
H.~N. Mhaskar.
\newblock Eignets for function approximation on manifolds.
\newblock {\em Applied and Computational Harmonic Analysis}, 29(1):63--87,
  2010.

\bibitem{heatkernframe}
H.~N. Mhaskar.
\newblock A generalized diffusion frame for parsimonious representation of
  functions on data defined manifolds.
\newblock {\em Neural Networks}, 24(4):345--359, 2011.

\bibitem{fasttour}
H.~N. Mhaskar and J.~Prestin.
\newblock Polynomial frames: a fast tour.
\newblock {\em Approximation Theory XI: Gatlinburg}, pages 101--132, 2004.

\bibitem{videomoeslund2012}
T.~B. Moeslund.
\newblock {\em Introduction to video and image processing: Building real
  systems and applications}.
\newblock Springer Science \& Business Media, 2012.

\bibitem{mousazadeh2015embedding}
S.~Mousazadeh and I.~Cohen.
\newblock Embedding and function extension on directed graph.
\newblock {\em Signal Processing}, 111:137--149, 2015.

\bibitem{nakatsukasa2010optimizing}
Y.~Nakatsukasa, Z.~Bai, and F.~Gygi.
\newblock Optimizing {H}alley's iteration for computing the matrix polar
  decomposition.
\newblock {\em SIAM Journal on Matrix Analysis and Applications},
  31(5):2700--2720, 2010.

\bibitem{nikolskii}
S.~M. Nikolskii.
\newblock {\em Approximation of functions of several variables and imbedding
  theorems}.
\newblock Springer Verlag, 1975.

\bibitem{roweis2000nonlinear}
S.~T. Roweis and L.~K. Saul.
\newblock Nonlinear dimensionality reduction by locally linear embedding.
\newblock {\em Science}, 290(5500):2323--2326, 2000.

\bibitem{rudinreal}
W.~Rudin.
\newblock {\em Real and complex analysis}.
\newblock Tata McGraw-Hill Education, 1987.

\bibitem{rudinfunctional}
W.~Rudin.
\newblock Functional analysis. international series in pure and applied
  mathematics, 1991.

\bibitem{sikora2004riesz}
A.~Sikora.
\newblock Riesz transform, {G}aussian bounds and the method of wave equation.
\newblock {\em Mathematische Zeitschrift}, 247(3):643--662, 2004.

\bibitem{singer}
A.~Singer.
\newblock From graph to manifold {L}aplacian: The convergence rate.
\newblock {\em Applied and Computational Harmonic Analysis}, 21(1):128--134,
  2006.

\bibitem{szego}
G~Szeg\"o.
\newblock Orthogonal polynomials.
\newblock In {\em Colloquium publications/American mathematical society},
  volume~23. Providence, 1975.

\bibitem{tenenbaum2000global}
J.~B. Tenenbaum, V.~De~Silva, and J.~C. Langford.
\newblock A global geometric framework for nonlinear dimensionality reduction.
\newblock {\em Science}, 290(5500):2319--2323, 2000.

\bibitem{bertozzicommunity}
Y.~van Gennip, B.~Hunter, R.~Ahn, P.~Elliott, K.~Luh, M.~Halvorson, S.~Reid,
  M.~Valasik, J.~Wo, G.~E. Tita, et~al.
\newblock Community detection using spectral clustering on sparse geosocial
  data.
\newblock {\em SIAM Journal on Applied Mathematics}, 73(1):67--83, 2013.

\bibitem{watson}
G.~N. Watson.
\newblock {\em A treatise on the theory of Bessel functions}.
\newblock Cambridge university press, 1995.

\bibitem{weinberger2005nonlinear}
K.~Q. Weinberger, B.~D. Packer, and L.~K. Saul.
\newblock Nonlinear dimensionality reduction by semidefinite programming and
  kernel matrix factorization.
\newblock In {\em Proceedings of the tenth international workshop on artificial
  intelligence and statistics}, pages 381--388. Citeseer, 2005.

\bibitem{zhang2004principal}
Z.-Y. Zhang and H.-Y. Zha.
\newblock Principal manifolds and nonlinear dimensionality reduction via
  tangent space alignment.
\newblock {\em Journal of Shanghai University (English Edition)},
  8(4):406--424, 2004.

\end{thebibliography}

\end{document}